
\magnification 1200
\input amstex
\documentstyle{amsppt}
\NoBlackBoxes
\NoRunningHeads

\hsize = 6.7 truein

\define\ggm{{G}/\Gamma}

\define\vre{\varepsilon}
\define\hs{homogeneous space}
\define\df{\overset\text{def}\to=}
\define\la{\Lambda_{A}}

\define\br{\Bbb R}
\define\bn{\Bbb N}
\define\bz{\Bbb Z}

\define\bq{\Bbb Q}
\define\bc{\Bbb C}

\define\ba{badly approximable}

\define\da{Diophantine approximation}
\define\spr{Sprind\v zuk}
\define\thus{\quad\Rightarrow\quad}
\define\ssm{\smallsetminus}
\define\lol{logarithm law}

\define\ca{\Cal A}
\define\ch{\Cal H}
\define\cb{\Cal B}
\define\tdf{tail distribution function}

\define\aeto{\underset{\mu\text{-a.e.}}\to\to}

\define\shn{_{{}{\ch,N}}}

\define\san{_{{}{\ca,N}}}

\define\sn{_{{}{N}}}

\define\w{\bold\Upsilon}
\define\td{\Phi_{\Delta}}
\define\un#1#2{\underset\text{#1}\to#2}

\define\a{\goth a}

\redefine\d{\goth d}
\define\g{\goth g}
\define\p{\goth p}
\define\k{\goth k}

\define\mr{M_{m,n}(\br)}
\define\amr{$A\in M_{m,n}(\br)$}
\define\mn{{m+n}}
\define\ve{\bold e}
\define\vx{\bold x}
\define\vp{\bold p}
\define\vq{\bold q}
\define\vv{\bold v}
\define\vz{\bold z}
\define\ph{\varphi}
\define\vt{\bold t}

\define\nz{\smallsetminus \{0\}}


\topmatter
\title logarithm laws for flows on homogeneous spaces\endtitle

\author { D.$\,$Y.~Kleinbock} \\ 
  { \rm 
   Rutgers University} 
 \vskip 5pt 
 and G.$\,$A.~Margulis \\ 
  { \rm 
   Yale University} \\ \ 
\endauthor 

    \address{ Dmitry Y. Kleinbock, Department of Mathematics, Rutgers
University, New Brunswick, NJ 08903}
  \endaddress

\email kleinboc\@math.rutgers.edu \endemail


    \address{ G.$\,$A. Margulis, Department of Mathematics, Yale University, 
   New Haven, CT 06520}\endaddress

\email margulis\@math.yale.edu
  \endemail

\abstract 
In this paper we generalize and sharpen D.~Sullivan's \lol\ for
geodesics by specifying conditions on a sequence of subsets $\{A_{t}\mid
{t}\in \bn\}$ of a homogeneous space $\ggm$ ($G$ a semisimple Lie group,
$\Gamma$ an irreducible lattice)  
and a sequence of elements $f_{t}$
of $G$ under
which $\#\{{t}\in\bn\mid f_{t}x\in A_{t}\}$ is infinite for
a.e.~$x\in \ggm$. 
The main tool is exponential decay of correlation coefficients of
smooth functions on $\ggm$.  Besides the general (higher rank) version
of Sullivan's result, as a consequence we obtain a new proof of the
classical Khinchin-Groshev theorem on simultaneous \da, and settle a
conjecture recently made by M.~Skriganov.

\endabstract

 \thanks The work of the first named author was supported in part by NSF
Grants DMS-9304580 and DMS-9704489, and that of the second named
author by NSF Grants DMS-9424613 and DMS-9800607. \endthanks        

\date 
March 25, 1999 \enddate

\endtopmatter

\document

\heading{\S 1. Introduction
}\endheading

\subhead{1.1}\endsubhead This work has been motivated by the
following two related results. The first one is the 
Khinchin-Groshev theorem, one of the cornerstones of metric theory of \da.
 We will
denote by $\mr$ the space 
 of real matrices with  $m$ rows and $n$ columns, 
and by   $\|\cdot\|$ the norm on $\br^k$, $k\in\bn$,  given by $\|\vx\| = \max_{1\le i \le k}|x_i|$.

\proclaim{Theorem \rm \cite{G}} Let $m$, $n$ be positive integers and
$\psi:[1,\infty)\mapsto (0,\infty)$ a 
non-increasing continuous function. 
 Then for  almost every
(resp.~almost no)  $A\in\mr$  there are infinitely many $\vq\in \bz^n$
such that  
$$
 \|A\vq + \vp\|^m   \le \psi(\|\vq\|^n) \quad \text{for some
}\vp\in\bz^m\,,  \tag 1.1
$$
provided the  integral
$
{\int_{1}^\infty {\psi(x)}\,dx}
$
diverges (resp.~converges). 
\endproclaim

 \subhead{1.2}\endsubhead The second motivation comes from  
 the paper \cite{Su} of D.~Sullivan. Let $\Bbb H^{{k} +
1}$ stand for the ${{k} + 1}$-dimensional real hyperbolic space with
curvature $-1$.  Take a 
 discrete group $\Gamma$  of hyperbolic isometries of $\Bbb H^{{k} +
1}$ such that $Y = \Bbb H^{{k} + 1}/\Gamma$ is not compact and has
finite volume. For $y\in Y$, denote by $S_y(Y)$ the set of unit
vectors tangent to $Y$ at $y$, and by $S(Y)$ the unit tangent bundle
$\{(y,\xi)\mid{y\in Y},\,{\xi\in S_y(Y)}\}$  of $Y$. Finally, for
$(y,\xi)\in S(Y)$   let $\gamma_t(y,\xi)$ be the geodesic on $Y$
through $y$ in the direction of $\xi$. 
The following theorem is essentially proved in \cite{Su} (see Remark (1) in \S 9):

\proclaim{Theorem} For $Y$ as above, fix $y_0\in Y$,  and  let
$\{r_{t}\mid {t}\in \bn\}$ be an arbitrary  sequence of real numbers. Then
for any $y\in 
Y$ and  almost 
every (resp.~almost no) $\xi\in S_y(Y)$ there  are infinitely many
${t}\in \bn$
such that 
$$
\text{\rm dist}\big(y_0,\gamma_{t}(y,\xi)\big) \ge  r_{t}\,,\tag 1.2
$$
  provided the series 
$
{\sum_{{t}=1}^{\infty} e^{-{k} r_{t}}}
$
diverges (resp.~converges). \endproclaim

 \subhead{1.3}\endsubhead  A choice 
$
r_{t} = \frac1{\varkappa} \log {t}\,,
$
 where ${\varkappa}$ is arbitrarily close to $k$, 
 yields the following statement, which has been referred to as the {\sl logarithm law for geodesics\/}: 

\proclaim{ Corollary} For $Y$  as above,  any $y\in Y$ and  almost all
$\xi\in S_y(Y)$, 
$$
\limsup_{t\to\infty}\frac{\text{\rm
dist}\big(y,\gamma_t(y,\xi)\big)}{\log t} = 1/{k} \,.\tag 1.3 
$$
\endproclaim

 \subhead{1.4}\endsubhead It seems natural to ask whether one can
generalize the statements of Theorem 1.2 and Corollary 1.3 to other
locally symmetric spaces of noncompact type. On the other hand,
Sullivan used a geometric proof of the case $m = n = 1$ of Theorem 1.1
to prove Theorem 1.2; thus one can ask whether there exists  a  connection
between the general case of  the 
Khinchin-Groshev theorem and some higher rank analogue of Sullivan's
result. 

In this paper we answer both questions in the affirmative. In particular,
the following generalization of Sullivan's results can be proved:

\proclaim{Theorem} For any noncompact irreducible\footnote{In fact the
theorem is true for reducible spaces as well, see
\S 10.2 for details.} 
locally symmetric space $Y$ of
noncompact type and finite volume there exists ${k}  = {k} (Y)
> 0$ such that the following holds: if $y_0\in Y$ and  $\{r_{t}\mid {t}\in
\bn\}$ is  an arbitrary  sequence of positive numbers, then   for  any
$y\in Y$ and almost 
every (resp.~almost no) $\xi\in S_y(Y)$ there  are infinitely many
${t}\in \bn$
such that {\rm (1.2)} is satisfied, 
  provided the series $\sum_{{t}=1}^{\infty} e^{-{k} r_{t}}$
diverges (resp.~converges). Consequently, {\rm (1.3)} holds for any $y\in Y$
and  almost all  $\xi\in S_y(Y)$. 
\endproclaim

The constant ${k}(Y)$ can be explicitly calculated in any given
special case; in fact, 
$
{k}(Y) =$ 
$\lim_{r\to\infty}{-\log\big(\text{vol}(A(r))\big)}/r\,,
$
where   
$$
A(r)\df \{y\in Y\mid\text{\rm dist}(y_0,y) \ge  r\}\,, \tag 1.4
$$
and ``vol'' stands for a Riemannian volume.  In
other words, the series $\sum_{{t}=1}^{\infty} e^{-{k} r_{t}}$ is, up to a constant, the sum of volumes
of sets $A(r_t)$. The latter sets can be viewed as a ``target
shrinking to $\infty$'' (cf.~\cite{HV}), and Theorems 1.2 and 1.4 say
that  if the shrinking is  slow enough
(read: the sum of the volumes is infinite), then almost all geodesics
approach infinity faster than the sets $A(r_t)$. 


This ``shrinking target'' phenomenon, being one of the main themes of
the present paper, deserves an additional 
discussion. Thus we have to make a terminological digression.   Let
$(X,\mu)$ be a probability space 
and let $F = \{f_t\mid t\in \bn\}$ be a sequence of $\mu$-preserving
transformations of $X$. Also let $\cb$ be a  family
of measurable 
subsets of $X$.

\example{1.5. Definition} Say that $\cb$ is a {\sl
Borel-Cantelli family for $F$\/}  if
for every sequence $\{A_t\mid t\in \bn\}$ of sets from $\cb$ one has
$$
\mu\big(\{x\in X\mid f_t(x)\in A_t\text{ for infinitely many
}t\in\bn\}\big)= \cases 0 \quad\text{ if } &\sum_{t = 1}^\infty
\mu(A_t) < \infty\\ 1\quad\text{ if }&\sum_{t = 1}^\infty \mu(A_t) =
\infty 
\endcases
$$
Note that the statement on top is always true in view of the classical
Borel-Cantelli Lemma, see \S 2.3.  An important special case is $F =
\{f^t\mid t\in \bn\}$ for a  measure-preserving transformation
$f:X\mapsto 
X$. We will say that  $\cb$  is {\sl Borel-Cantelli for  $f$\/} if it is
Borel-Cantelli  for $F$ as above.  
\endexample

It is easy to see that  $f:X\mapsto
X$ is ergodic (resp.~weakly mixing\footnote{This characterization of
weak mixing was pointed out to us by Y.~Guivarc'h and A.~Raugi; see
also \cite{CK}.}) iff
every one-element  (resp.~finite) family of sets of positive measure
is Borel-Cantelli for $f$.  
On the other hand, if $(X,\mu)$ is nontrivial,
then for any sequence of transformations $F = \{f_t\}$ one can
construct a family  (say, $A_t = f_t(A)$ with $0< \mu(A) <
1$) which  is not Borel-Cantelli for $F$. Therefore in order to describe Borel-Cantelli
families of sets for a particular sequence of maps, it is natural to
specialize and impose certain regularity restrictions on the sets
considered.  

An important example is given in the paper \cite{P} of W.~Philipp:
there $X = [0,1]$, $f$ is an expanding map of $X$ given by either
$x\mapsto \{\theta x \}$, $\theta > 1$, or by $x\mapsto \{\frac1x\}$
($\{\cdot\}$ stands for the fractional part), and it is proved that
the family  of all intervals is Borel-Cantelli for $f$. This means
that one can take any $x_0\in [0,1]$ and consider a ``target shrinking
to $x_0$'', i.e.~a sequence  $(x_0 - r_t, x_0 +
r_t)$. Then almost all orbits $\{f^tx\}$ get into infinitely many such
intervals whenever $r_t$ decays  slowly enough. This can be
thought of as a quantitative strengthening of density of almost all orbits
(cf.~the paper \cite{Bos} for a similar approach to the
rate of recurrence). 

We postpone further discussion of this general set-up until 
 \S 10.2, and concentrate on ``targets shrinking to infinity'' in
noncompact  spaces. Our goal is to state a result which will
imply both Theorem 1.4 and Theorem 1.1.  For $Y$ as in Theorem 1.4,
let $G$ be the connected 
component of the identity in the isometry 
group of the universal cover of $Y$. Then $G$ is a connected semisimple Lie
group  without compact factors, and the space $Y$ can be identified
with $K\backslash \ggm$, where 
$\Gamma$  is an irreducible
lattice in $G$ and $K$ is a maximal compact subgroup of $G$. Instead of
working with $Y$, we choose
the \hs\ $X = \ggm$ as our main object of investigation. Fix a Cartan
subalgebra  $\a$ of the Lie algebra of $G$.
It is known 
\cite{Ma} that the geodesic flow on  the unit tangent bundle $S(Y)$
of $Y$ can be realized via  action of 
one-parameter subgroups of the form $\{\exp(t\vz)\}$, with $\vz\in\a$,
on the space $X$ (see \S 6 for
details).   In what
follows, we will choose a  maximal compact subgroup $K$ of $G$,
 endow $X$ with a Riemannian metric by fixing  a right invariant
Riemannian metric 
on $G$ bi-invariant with respect to
$K$, and let $\mu$ be the  normalized Haar measure  on $X$.

Recall that the ``neighborhoods of $\infty$'' of Theorem 1.4 are the
complements $A(r)$, see (1.4), of balls in $Y$, and it follows from
that  theorem  
 that the family $\{A(r)\mid r > 0\}$ is Borel-Cantelli for the time-one map of
the geodesic flow. To 
describe sequences of sets ``shrinking to infinity'' in $X$, we 
will   
replace  the distance function $\text{\rm dist}(y_0,\cdot)$ by a
function $\Delta$ on $X$ satisfying certain properties, and consider
the family  
 $$
{\cb}(\Delta)\df\big\{\{x \in X\mid \Delta(x) \ge r\}\mid
r\in\br\big\}
$$
of super-level sets of $\Delta$. To specify the class of functions
$\Delta$ that we will work with,  we
introduce the following 

\example{1.6. Definition}  For a function $\Delta$ on $X$,  define the
{\sl \tdf\/}  
$\td$ of $\Delta$ by  
$$
\td(z) \df {\mu\big(\{x\mid \Delta(x)\ge z\}\big)}\,.
$$ 
Now   say that 
$\Delta$ is {\sl DL\/} (an abbreviation for
``distance-like'') if it is  uniformly continuous, and
$\td$ does not decrease very fast, more precisely, if
$$
\exists\, c,\delta > 0
\text{ such that }\td(z + \delta) \ge c\cdot\td(z)\ \forall z\ge 0\,.\tag DL
$$
For $k > 0$, we will also say that  $\Delta$ is {\sl $k$-DL\/}
if it is  uniformly continuous and in addition 
$$
\exists\,C_1,C_2>0\text{ such that } C_1e^{-kz} \le \td(z)\le C_2e^{-kz} \quad\forall\,z\in\br\,.\tag $k$-DL
$$
\endexample

It is clear  that ($k$-DL) implies (DL). Note  that DL functions on
$X$ exist only when $X$ is not compact (see
\S 4.3).
The most important example (\S 5) is the distance function on $X$. 
%
Thus the following theorem can be viewed as a generalization of
Theorem 1.4:


\proclaim{1.7. Theorem} Let $G$ be a connected semisimple Lie
group  without compact factors,  $\Gamma$ an irreducible
lattice in $G$, 
$\a$ a Cartan
subalgebra of the Lie algebra of $G$, $\vz$ a nonzero element of
$\a$. Then:

{\rm (a)}  if $\Delta$ is a DL function  on $X = \ggm$,   the family
${\cb}(\Delta)$ is  Borel-Cantelli for $\exp(\vz)$; 

{\rm (b)} if  $\Delta$ is 
 $k$-DL for some $k > 0$, then   for almost all $x\in X$ one has
$$
\limsup_{t \to +\infty}
\frac{\Delta\big(\exp({t\vz})x\big)}{\log t} = 1/{k} \,.\tag 1.5
$$
\endproclaim

In particular,  (1.3) can be derived from (1.5) by taking $G =
SO_{{k} + 1,1}(\br)$ 
and $\Delta(x) =  \text{dist}(x_0,x)$ for fixed  $x_0
\in  \ggm$.    

\subhead{1.8}\endsubhead In fact, it
is possible to derive a version of Theorem 1.7 for actions of
multi-parameter subgroups of $G$. More generally,  we will consider
actions  
of arbitrary countable sequences $\{f_{t}\mid 
{t}\in \bn\}$  of elements of $G$.
To specify a class of sequences good for our purposes, 
denote by $\|g\|$ the
distance between $g\in G$ and the identity element of $G$, and say
that  a sequence $\{f_{t}\}$ is {\sl ED\/} (an abbreviation for
``exponentially divergent'') if 
$$
\sup_{t\in\bn}\sum_{s=1}^\infty e^{-\beta\|f_sf_t^{-1}\|} < \infty \quad
\forall\,\beta > 0\,.
\tag ED 
$$
In this setting we state the following general result:


\proclaim{Theorem} 
 For $G$ and $\Gamma$  as in Theorem 1.7, let $F = \{f_{t}\mid
{t}\in \bn\}$ be an  
ED  sequence  of elements of $G$ and  $\Delta$ a
DL function on $\ggm$.  Then   the family
${\cb}(\Delta)$ is Borel-Cantelli  for $F$. \endproclaim

\subhead{1.9}\endsubhead Clearly Theorem 1.7 is
 a special case of the above theorem:  it is easy to check 
(see \S
4.4) 
that  the
sequence $f_{t} 
= \exp({{t}\vz})$, with $\vz\in\a\nz$,  satisfies (ED). 
More
generally, the following multi-parameter generalization of Theorem 1.7 can
be derived from Theorem 1.8: 

\proclaim{Theorem} For  $G$, $\Gamma$, $X$  and   $\a$ as in Theorem
1.7,

{\rm (a)} if $\Delta$  is a DL function  on $X$,
and ${t}\mapsto \vz_{t}$ is a map from  $\bn$ to 
$\a$ such that 
$$
\inf_{{t}_1\ne {t}_2}
\|\vz_{{t}_1}-\vz_{{t}_2}\| > 0\,,\tag 1.6
$$ 
 then    the family
${\cb}(\Delta)$ is Borel-Cantelli for $\{\exp(\vz_{t})\mid {t}\in
\bn\}$; 

{\rm (b)} if  $\Delta$ is 
 $k$-DL for some $k > 0$, and  $\d_{\sssize +}$ is  a  nonempty open
cone in a $d$-dimensional 
subalgebra $\d$  of $\a$ 
($1\le d \le \text{\rm rank}_\br(G)$), then   
for almost all $x\in X$ one
has 
$$
\limsup_{\vz \in \d_{\sssize +},\,\vz\to\infty}
\frac{\Delta\big(\exp({\vz})x\big)}{\log \|\vz\|} = d/{k} \,.\tag 1.7
$$
 \endproclaim

\subhead{1.10}\endsubhead From the above theorem one can get \lol s
for flats in locally 
symmetric spaces. Let the space $Y$ be as in Theorem 1.4. As usual, by
a $d$-dimensional flat in   $Y$ ($1\le d \le
\text{\rm 
rank}(Y)$) we mean the image of $\br^d$ under a
locally isometric embedding into $Y$. For $y\in Y$,  denote by
$S_y^d(Y)$ the set of orthonormal $d$-tuples  of vectors  $\xi_i\in
S_y(Y)$ which form a basis for a tangent space to a flat passing
through $y$. The set $S_y^d(Y)$ is a real algebraic variety coming
with the natural measure class, which makes it possible to talk about
``almost all flats passing
through $y$''. If $\vec \xi = 
(\xi_1,\dots,\xi_d)\in S_y^d(Y)$, we will denote by $\vt =
(t_1,\dots,t_d) \mapsto  \gamma_\vt(y,\vec\xi)$ the embedding
specified by $\vec \xi$, that is, we let $\gamma_\vt(y,\vec\xi) \df \exp_y(\sum_i t_i
\xi_i)$ (a multi-dimensional analog of the geodesic in the direction
of a single vector $\xi \in S_y(Y)$).

\proclaim{Theorem} Let  $Y$, $y_0$ and ${k}  = {k} (Y)$
be as in Theorem 1.4. Take  $1\le d \le
\text{\rm 
rank}(Y) $ and  a nonempty  open cone
$\d_{\sssize +}\subset \br^d$, and let $\vt\mapsto r_{\vt}$, $\vt\in
\d_{\sssize +}\cap \bz^d$, be a real-valued function. 
Then for any $y\in Y$ and almost 
every (resp.~almost no) $\vec\xi\in S_y^d(Y)$ there  are infinitely many
${\vt}\in \d_{\sssize +}\cap \bz^d$
such that 
$
\text{\rm dist}\big(y_0,\gamma_{\vt}(y,\vec\xi)\big) \ge  r_{\vt}\,,
$ 
  provided the series 
$
{\sum_{\vt\in
\d_{\sssize +}\cap \bz^d} e^{-{k} r_{\vt}}}
$
diverges (resp.~converges). Consequently, for any
$y\in Y$ and
almost all $\vec \xi\in S_y^d(Y)$ 
one has  
$$
\limsup_{\vt \in \d_{\sssize +},\,\vt\to\infty}\frac{\text{\rm
dist}\big(y,\gamma_\vt(y,\vec \xi)\big)}{\log \|\vt\|} =
d/{k} \,.\tag 1.8
$$
\endproclaim


\subhead{1.11}\endsubhead Another class of applications of Theorems 1.7
and 1.9 
is given by  a modification  of 
S.G.~Dani's \cite{D, \S 2} correspondence between \da\ of  systems   of
$m$ linear forms in $n$ variables 
and flows on 
the space of lattices in $\br^k$, where $k = \mn$. Namely, consider $G =
SL_{k}(\br)$, $\Gamma = SL_{k}(\bz)$, and the function $\Delta$ on
the space  $\ggm$ of unimodular lattices in $\br^k$  defined by
$$
\Delta(\Lambda)  \df
\max_{\vv\in\Lambda\nz}\log\big(\tfrac1{\|\vv\|}\big)\,.\tag 1.9
$$
Denote also by $f_{t}$ the
element of $G$ of the form 
$$
f_{t} = \text{\rm diag}(\underbrace{e^{{t}/m},\dots,e^{{t}/m}}_{\text{$m$
times}},\underbrace{e^{-{t}/n},\dots,e^{-{t}/n}}_{\text{$n$
times}})\,.\tag 1.10
$$ 

We will show in \S 8  that  Theorem 1.1 follows from the fact that
the family   
${\cb}(\Delta)$ is Borel-Cantelli for $f_{1}$. 
Using similar technique, one can also prove a result that was, in somewhat
weaker form, conjectured by M.~Skriganov in \cite{Sk}:

\proclaim{Theorem} Let $\psi:[1,\infty)\mapsto (0,\infty)$ be a 
non-increasing continuous function and $k$ an integer greater than $1$.
Then for  almost every (resp.~almost no) unimodular lattice 
$\Lambda$ in $\br^k$  there are infinitely many $\vv \in \Lambda$ such that 
$$
\Pi(\vv) \le \|\vv\|\cdot\psi(\|\vv\|)   \tag 1.11
$$
(here and hereafter we use the notation $\Pi(\vv) \df \prod_{i = 1}^k
|v_i|$ for  $\vv = (v_1,\dots,v_k)\in \br^k$), provided the integral
$
{\int_{1}^\infty {(\log x)^{k -2}\psi(x)}\,dx} 
$
diverges (resp.~converges). 
\endproclaim

In \S 9 we will explain why the above statement can be thought of as a
higher-dimensional multiplicative
generalization of Khinchin's Theorem, and how one can derive it from
Theorem 1.9  by considering the
action of the whole Cartan subgroup of $SL_k(\br)$ 
on the space  $SL_{k}(\br)/SL_{k}(\bz)$.

The paper is organized as follows. In \S  2  we work in
a general setting of a probability  space $(X,\mu)$ and a  sequence of
nonnegative measurable functions $\ch = \{h_{t}\mid {t}\in\bn\}$ on
$X$, and,
following V.~\spr, 
write down a condition (Lemma 2.6)  
which guarantees that for almost every $x\in X$ the sum  $\sum_{{t} =
1}^\infty h_{t}(x)$ is infinite.   Then we throw in a
measure preserving action of $F = \{f_{t} \mid {t}\in \bn\}$ and apply the
aforementioned results to the {\sl twisted sequence\/} $\ch^F \df\
\{f_{t}^{-1}h_{t}\}$.   

In \S 3 we restrict ourselves to flows on $\ggm$ and prove the
following 

\proclaim{1.12. Theorem}  Let $G$ be a connected semisimple center-free 
Lie 
group  without compact factors,  $\Gamma$ an irreducible
lattice in $G$,
and let $\rho_0$  stand for the regular
representation of $G$ on the subspace  of $L^2(\ggm)$
orthogonal to 
constant functions. Assume in addition that $\ggm$ is not
compact. Then  the restriction   of $\rho_0$ to any simple 
factor of $G$
is isolated (in the Fell topology) from the
trivial representation. 
\endproclaim

The latter condition is known (cf.~\cite{KM, \S 2.4}) to  guarantee
exponential decay of correlation coefficients of smooth functions on
$\ggm$, see Corollary 3.5. In the next section we use  the fact that
$\Delta$ is DL 
to 
approximate characteristic functions of the sets $\{x\in\ggm\mid
\Delta(x)\ge r_{t}\}$ by smooth functions $h_t$. A
quantitative strengthening of Theorem
1.8 is then proved by deriving \spr's
condition  from estimates on decay of  correlation
coefficients of functions $h_t$.  Theorem 1.9  (hence 
1.7 as well)   is also 
proved in \S 4.   After that we  describe applications to 
geodesics and flats in locally symmetric spaces (Theorems 1.4 and 1.10)
and \da\ (Theorems 1.1 and 1.11).

 \heading{\S2. Borel-Cantelli-type results}\endheading 

\subhead{2.1}\endsubhead  Let $(X,\mu)$ be a probability space. We
will use  notation $\mu(h)=\int_X h\, d\mu$ for an
integrable function $h$ on $X$. Let us consider sequences $\ch =
\{h_{t}\mid {t}\in \bn\}$  of 
nonnegative integrable \footnote{Throughout the sequel
all the functions $h_{t}$ will be assumed measurable, integrable,
a.e.~nonnegative and nonzero on a set of positive measure.} functions
on $X$, and, for $N = 1,\dots,\infty$,  denote 
$$
S\shn(x) \df \sum_{{t} = 1}^N h_{t}(x)\quad\text{and}\quad E\shn \df \dsize\sum_{{t} = 1}^N \mu(h_{t})  =  \mu(S\shn) \,;
$$
this notation will be used throughout the paper. 
We will omit the index $\ch$ when it is clear from the context.
A special case of such a sequence is given by 
characteristic functions $h_{t} = 1_{ A_{t}}$, where $\ca = \{A_{t}
\mid {t}\in \bn\}$ is a sequence of measurable subsets of $X$. In this
case we will put the index $\ca$ in place of $\ch$ in the above
notation. 
We will  say that  a sequence $\ch$ (resp.~$\ca$) of functions
(resp.~sets) is {\sl summable\/} if  $E_{{ \ch},\infty}$ (resp.~$E_{{
\ca},\infty}$) is finite, and {\sl nonsummable\/} otherwise.

\example{Main example} If $\Delta$ is any 
function on
$X$ and $\{r_{t}\mid {t}\in \bn\}$
 a  sequence of real numbers, one can consider the sequence of
super-level sets 
$\{x \mid \Delta(x) \ge r_{t}\}$ of $\Delta$;
their measures are equal to $\td{(}r_{t}{)}$,  where $\td$ is the
\tdf\ (see \S 1.6) of $\Delta$. \endexample 

\example{2.2. Another main example}  Let $F = \{f_{t}\mid {t}\in \bn\}$ be a sequence of $\mu$-preserving  transformations of $X$. Then given any sequence $\ch = \{h_{t}\mid {t}\in \bn \}$ of functions on $X$ or a sequence $\ca = \{A_{t}\mid {t}\in \bn \}$ of subsets of $X$, one can consider {\sl twisted\/} sequences 
$$
\ch^F \df \{f_{t}^{-1}h_{t}\mid {t}\in \bn \}\quad\text{and}\quad \ca^F \df \{f_{t}^{-1}A_{t}\mid {t}\in \bn \}\,.
$$
By $F$-invariance of $\mu$, $E_{\ch^F, N}$ is the same as $E_{\ch,
N}$ for any $N\in\bn$; in particular, the twisted sequence is summable
if and only if the original one is. \endexample

\subhead{2.3}\endsubhead Given a sequence $\ca = \{A_{t} \mid {t}\in \bn
\}$ and a $\mu$-generic point $x\in X$, one may want to look at the
asymptotics of $S\san(x) = \#\{1\le {t}\le N\mid x\in A_{t}\}$ in comparison with the sum $E\san$ of measures of the sets $A_{t}$, ${1\le {t}\le N}$, as $N\to\infty$. This is for example the subject of the classical Borel-Cantelli Lemma. In general, for a  sequence $\ch$ of functions on $X$, it is very easy to estimate the ratio of $S\shn(x)$ and $E\shn$ from above as follows:

\proclaim{Lemma \rm  (cf.~\cite{KS, part (i) of the Theorem})} Let $(X,\mu)$ be a probability space, $\ch$ a sequence of functions on $X$. Then 
$$
\liminf_{N\to\infty} \frac{S\shn(x)} {E\shn} < \infty\quad\text{for }\mu\text{-a.e. } x\in X\,.
$$  
In particular, if  $\ch$  is summable, $S_{{ \ch},\infty}$ is  finite almost everywhere. 
 \endproclaim
 
\demo{Proof} By the Fatou Lemma,
$
\mu\Big(\liminf_{N\to\infty} \frac{S\shn} {E\shn} \Big) \le
\liminf_{N\to\infty}\mu\big(\frac{S\shn} {E\shn}\big) = 1 
$. \qed
\enddemo

One immediately recognizes the last assertion as  the conclusion of
the easy part of the classical Borel-Cantelli Lemma. It takes care of 
the convergence cases in all the Khinchin-type theorems stated in  the
introduction, 
as well as of the
 upper estimates for the limits in \lol s (1.3), (1.5), (1.7) and
(1.8).

\subhead{2.4}\endsubhead The corollary below will make the connection
with \lol s 
more transparent. We need the  following terminology: a real-valued function 
$r(\cdot)$ will be called {\sl quasi-increasing\/} if there exists a
constant $C$ such that 
$$
r(t_2) > r(t_1) - C\text{ whenever }t_1 \le t_2 <
t_1 + 1\,.\tag 2.1
$$ 

\proclaim{Corollary} Let $X$ be a metric space, $\mu$  a
probability measure on $X$, $d,k\in \bn$, $\d_{\sssize +}\subset \br^d$ a nonempty open
cone, $\vz\mapsto f_\vz$ a
continuous\footnote{Here by the distance between two maps
$f_1,f_2:X\mapsto X$ we mean $\sup_{x\in X}\text{dist}_X\big(f_1(x) -
f_2(x)\big)$.} homomorphism from $\d_{\sssize +}$ to the semigroup of
all  self-maps of $X$, $\Delta$ a $k$-DL 
function on $X$. For some  $t_0\in\br$, let   $r:[t_0,\infty)\mapsto
\br$ be a quasi-increasing 
 function  such that  the integral
$$ 
\int_{t_0}^\infty t^{d-1} e^{-kr(t)}\,dt \tag 2.2
$$
converges. Then for $\mu$-almost all $x\in X$ one has 
$
 \Delta\big(f_\vz (x)\big) < r(\|\vz\|)
$ 
whenever
$\vz\in\d_{\sssize +}$ is far  enough from $0$. Consequently,
$$
\limsup_{\vz \in \d_{\sssize +},\,\vz\to\infty}
\frac{\Delta\big(f_{\vz}(x)\big)}{\log \|\vz\|} \un{$\mu$-a.e.}{\le}
d/{k} \,.\tag 2.3
$$
\endproclaim
 

\demo{Proof} Choose a lattice  $\Sigma$ in $\br^d$;
from (2.1) and
the convergence of (2.2) it
follows that  the series
$$
\sum_{\vz\in\d_{\sssize +} \cap\Sigma,\,\|\vz\| \ge t_0}e^{-k r({\|\vz\|})} \tag 2.4
$$
converges. Clearly for  any $\vz\in\d_{\sssize +}$ far enough from $0$
one can find $\vz'\in \d_{\sssize +} \cap \Sigma$ such that 
$$
\|\vz\| - 1\le \|\vz'\|\le  \|\vz\|\,,\tag 2.5
$$
 and
$\|\vz' - \vz\|$ is less than some uniform constant $C_1$. Since 
the correspondence $\vz\mapsto f_\vz$ is continuous, for some $C_2$
one then has 
$\sup_{x\in X}\text{dist}\big(f_\vz(x),f_{\vz'}(x)\big) < C_2$;
further, from the uniform continuity of $\Delta$ it follows that for
some $C_3$ one has
$$\sup_{x\in X}|\Delta\big(f_\vz(x)\big) - \Delta\big( f_{\vz'}(x)\big)|
< C_3\,.\tag 2.6
$$ 

Now consider the sequence of sets
$
\ca \df \big\{ \{x\in X\mid \Delta(x) \ge r(\|\vz\|) - C - C_3\}
\bigm| \vz\in\d_{\sssize +} \cap 
\Sigma\big\}$, with 
$C$ as in (2.1), 
and $F = \{ f_{\vz}\mid \vz\in\d_{\sssize +} \cap \Sigma\}$. It
follows from the 
convergence of (2.4) and
$\Delta$ being DL that $\ca$ is summable. 
Applying Lemma 2.3 to the twisted sequence $\ca^F$, one
concludes that for almost all $x$ one has   $ \Delta\big(f_{\vz'}
(x)\big) < r(\|\vz'\|)-  C - C_3$ 
for $\vz'\in\d_{\sssize +}\cap \Sigma$ with large enough
$\|\vz'\|$. In view of (2.1), 
(2.5) and (2.6), this implies that for
almost all $x$  one has $
 \Delta\big(f_\vz (x)\big) < r(\|\vz\|)
$ 
  for all $\vz\in\d_{\sssize +}$ with large enough $\|\vz\|$. The
second part of the corollary is obtained by taking $r(t) = \frac
d{\varkappa} \log {t}$ with $\varkappa < k$. The integral (2.2)
obviously converges, therefore for almost all $x$ one has
$
\dsize{\frac{\Delta\big(f_{\vz}(x)\big)}{\log \|\vz\|} <  \frac
d{\varkappa} }
$
 whenever
$\vz\in\d_{\sssize +}$ is far  enough from $0$, and (2.3) follows. 
\qed
\enddemo

\example{2.5. Example} Take $X = S(Y)$ as in \S 1.2, $\mu$ the
Liouville measure on $S(Y)$, fix $y_0\in Y$ and  let
$\Delta\big((y,\xi)\big) = \text{dist}(y_0,y)$. As mentioned in 
 \cite{Su, \S 9}, $\Delta$ is $k$-DL. From the above corollary (with
$d = 1$ and $\d_{\sssize +} = \br_{\sssize +}$) one concludes  that 
$
\dsize{\limsup_{t\to\infty}\frac{\text{\rm
dist}\big(y_0,\gamma_t(y,\xi)\big)}{\log t} }
$
  as
$t\to\infty$ is not greater than $1/{k} $. To derive the upper
estimate for the limit in Corollary 1.3 from the above statement, it
suffices to observe that for any two points $y_1,\,y_2$ of $Y$:

\roster 
\item"{$\bullet$}" the functions dist$(y_1,\cdot)$ and
dist$(y_2,\cdot)$ differ by at most dist$(y_1,y_2)$, and

\item"{$\bullet$}" for any geodesic ray $\gamma$ starting from $y_1$ there
is a geodesic  ray starting from $y_2$ which stays at a bounded distance
from $\gamma$. \endroster 
%
\endexample


\subhead{2.6}\endsubhead Let $F$ be a sequence
of $\mu$-preserving 
transformations of $X$ and  ${\cb}$ a family of measurable
subsets of $X$. From Lemma 2.3 it is clear that ${\cb}$ is Borel-Cantelli for $F$
iff for any    nonsummable
sequence $\ca$ of sets from $\cb$ one has
 $S_{{ \ca^F},\infty} = \infty$ for almost all $x\in
X$. Therefore we are led to 
studying asymptotical lower  estimates for ${S\shn} /{E\shn}$, with
$\ch$ as in \S 2.1. 

 One can easily find many examples of sequences $\ch$ for which the
above ratio almost surely 
tends to zero 
as $N\to\infty$. It is also well known (see \cite{Sp, p.~317} for a
historical overview) that  
 the estimates we are after 
 follow from certain conditions on second moments of the functions $h_{t}$. 
We will  employ a lemma which was abstracted by V.~\spr\
 from the works of W.~Schmidt (see also \cite{P}
for a related result).

\proclaim{Lemma \rm(\cite{Spr, 
Chapter I,  Lemma 10})} For  a  sequence   $\ch = \{h_t\mid t\in\bn\}$
of functions on 
$X$, assume that  
$$
\mu(h_{t}) \le 1\quad\text{for all }t\in \bn \tag 2.7
$$
and 
$$ 
\exists\,C > 0\text{ such that }\int_X\Big(\sum_{{t} = M}^{N}
h_{t}(x) - \sum_{{t} = M}^{N} \mu(h_{t}) \Big)^2\,d\mu  \le C\cdot
\sum_{{t} = M}^{N} \mu(h_{t})  \quad \forall\,N> M \ge 1\,.\tag SP
$$ 
Then  for any positive $\vre$ one has, as $N\to\infty$,  
$$
S\shn(x) = {E\shn} + O(E\shn^{1/2} \log^{3/2 + \vre} E\shn)
\tag 2.8 
$$  for
$\mu$-a.e.~$x\in X$; in particular, 
$\dfrac{S\shn(x)}{E\shn} \underset{\text{a.e.}}\to\to 1\text{ as
}N\to\infty
$ whenever $\ch$ is  nonsummable. 
\endproclaim

\comment 

Note that this is is a quasi-independent analogue of Theorem 6.6 from
\cite{Du}; our proof  essentially follows the lines of Durrett's
argument.  

\demo{Proof} We will drop the subscript $\ch$ and introduce the
notation $\Psi\sn = 
\frac{S\sn}{E\sn}$. Note that (SP) can be written in the form 
$
\text{Var}(\Psi\sn) \le  C/{E\sn}
$. Therefore,  by Chebyshev's inequality, for any $\vre > 0$ and  $N\in\bn$
$$
\mu\big(\{x\in X \bigm| |\Psi\sn(x) - 1| > \vre \}\big) \le
\frac{\text{Var}(\Psi\sn)}{\vre^2} \le  \frac
C{E\sn\vre^2} \,.\tag 2.7
$$
Define an increasing sequence $\{N_{t}\}\subset \bn$, ${t}\in\bn$, by
$$
N_{t} = \min\{N\mid E\sn \ge {t}^{2}\}\,.\tag 2.8
$$
It follows from (2.7) that for any $\vre > 0$
the family $\{x\in X \mid |{\Psi{_{{}{N_{t}}}}} -
1| > \vre\}$ is summable, so, by Lemma 2.3,
$
\#\{{t}\in\bn\mid | {\Psi{_{{}{N_{t}}}}}- 1| > \vre\}$ is finite for $\mu\text{-a.e. } x\in X$.
Since $\vre$ is arbitrarily small, it follows that $\Psi_{{}{N_{t}}}\aeto 1$ as ${t}\to\infty$.

Now for any $N\in\bn$ one finds ${t}$ such that $N_{t}\le N <
N_{{t}+1}$ and, using  (2.8) and the fact that all $\alpha_{t}$ are
majorated by some number $L$, writes  
$$
\frac{{t}^2}{({t}+1)^2 + L}{\Psi{_{{}{N_{t}}}}} \le \frac{E{_{{}{N_{t}}}}}{E{_{{}{N_{{t}+1}}}}}{\Psi{_{{}{N_{t}}}}} \le \Psi\sn \le
\frac{E{_{{}{N_{{t}+1}}}}}{E{_{{}{N_{t}}}}}{\Psi{_{{}{N_{{t}+1}}}}}
\le \frac{({t}+1)^2 + L}{{t}^2}{\Psi{_{{}{N_{{t}+1}}}}}\,.
$$
This shows that  $\mu$-almost everywhere one has $\Psi\sn \to 1$ as
$N\to\infty$. \qed\enddemo 

\endcomment

\example{2.7. Remark}  Note that the left hand side of  (SP) is equal to 
$$
 \int_X\biggl(\sum_{{t} = M}^{N} h_{t}\biggr)^2d\mu -  \biggl(\dsize\sum_{{t} = M}^N \mu(h_{t})\biggr)^2 
 = \sum_{{s,t}  = M}^{N} \big(\mu(h_sh_t) - \mu(h_{s})\mu( h_{t})
\big)\,.\tag 2.9 
$$
This shows that (SP) will hold provided the correlation coefficients
$|\mu(h_sh_t) -  \mu(h_{s})\mu( h_{t})|$ become small for large values of
$|s-t|$. Our plan is to apply Lemma 2.6 to the twisted
sequences $\Cal H^F$, where $F$ is as in Theorem 1.8 and $\Cal H$
consists of smooth functions on $\ggm$. The exponential  decay of
correlations under the 
$F$-action, the  main result of the next section, will be enough to
guarantee  (SP).
\endexample

\subhead{2.8}\endsubhead We close the section with a partial converse
to Corollary 2.4, which we will use later for the derivation of
\lol s.

\proclaim{Lemma} Let $X$, $\mu$, $d$, $k$, $\d_{\sssize +}$, $\vz\mapsto f_\vz$,
$\Delta$ and $t_0$ be as in Corollary 2.4, and let  $r:[t_0,\infty)\mapsto
\br$ be  a quasi-increasing 
 function   such that  the integral
{\rm (2.2)} diverges.  Assume that there exists a lattice  $\Sigma$ in
$\br^d$ such that the family $\cb(\Delta)$ of super-level sets of
$\Delta$ 
is Borel-Cantelli for $F \df\{f_\vz 
\mid\vz\in\d_{\sssize +} \cap \Sigma\}$.  
 Then  for $\mu$-almost all $x\in X$ there exist
$\vz\in\d_{\sssize +}$ arbitrarily far from $0$ such that 
$
 \Delta\big(f_\vz (x)\big) \ge  r(\|\vz\|)\,.
$ 
Consequently,
$
\dsize{\limsup_{\vz \in \d_{\sssize +},\,\vz\to\infty}
\frac{\Delta\big(f_{\vz}(x)\big)}{\log \|\vz\|} \un{$\mu$-a.e.}{\ge}
d/{k} \,.}
$
\endproclaim
 

\demo{Proof} From (2.1) and
the divergence of (2.2) it
follows that  the series (2.4) is  divergent. In view of $\Delta$ being
$k$-DL and by definition of $\cb(\Delta)$ 
being  Borel-Cantelli for $F$, one gets $
 \Delta\big(f_\vz (x)\big) \ge  r(\|\vz\|)$  almost surely for infinitely
many $\vz\in\d_{\sssize +} \cap 
\Sigma$, hence the first part of the lemma. The 
second part is immediate by taking  $r(t) = \frac
d{k} \log {t}$.
\qed
\enddemo

\heading{\S3. Isolation properties of representations and correlation
decay}\endheading


\subhead{3.1}\endsubhead Let $G$ be a locally compact second countable
group. Recall that the Fell topology
on the set of (equivalence classes of) unitary representations $\rho$ of $G$
 in  separable Hilbert spaces $V$ is defined so that the sets $\{ \rho
\bigm| \|\rho(g)v - v\| < \vre \| 
v\|\ \forall g\in K \ \forall v\in V\}$, where $\vre > 0$ and $K$ runs
through all
compact subsets of $G$, constitute a basis of open neighborhoods of the
trivial representation $I_G$ of $G$. (See the Appendix  and \cite{M,
Chapter III} for more detail.)  If $(X,\mu)$ is a probability space
and $(g,x)\mapsto gx$ a 
 $\mu$-preserving action of $G$ on $X$, we will denote by $L^2_0(X,\mu)$
the subspace  of $L^2(X,\mu)$ 
orthogonal to 
constant functions. 
Our proof of Theorem 1.12 will use the following
result, communicated by A.~Furman and Y.~Shalom, which will allow us
to pass from a space to its finite covering:

\proclaim{Lemma} Let  $(X_1,\mu_1)$  and $(X_2,\mu_2)$ be probability
spaces, $G$ a locally compact
second countable  group
acting ergodically on both,  and let  $\pi:X_1\mapsto X_2$ be a
surjective measurable $G$-equivariant map such that for
some positive $c<1$ one has 
$$
c \mu_1(A) \le \mu_2\big(\pi(A)\big) \le \frac1c
\mu_1(A)\quad\text{for any }A\subset X_1\,.\tag 3.1
$$
  Denote
by  $\rho_{i,0}$ 
the regular
representation of $G$ on $L^2_0(X_i,\mu_i)$ 
($i = 1,2$).
Then $\rho_{1,0}$ is isolated from $I_G$ iff so is
$\rho_{2,0}$.  
  \endproclaim

The proof of Furman and Shalom is based on the connection between
$\rho_0$ being close to  $I_G$ and existence of nontrivial
$G$-invariant means on $L^\infty(X,\mu)$ \cite{FS, Theorem 1.8}. In
the Appendix we  give a more transparent proof, based on the notion of
asymptotically invariant sequences of subsets of $X$. The argument
goes back to J.~Rosenblatt \cite{Ro} and K.~Schmidt \cite{S} and runs
more or less in parallel to the proof given in  \cite{FS}.

\subhead{3.2}\endsubhead  Let now $G$ be a connected semisimple center-free Lie
group without compact factors,  $\Gamma$ an irreducible
lattice in $G$, $\mu$ the normalized Haar measure on the \hs\ $\ggm$.
  It is known (see \cite
{B, Lemma 3}) that  the regular 
representation $\rho_0$  of $G$  on  $L^2_0(\ggm,\mu)$ is isolated from
$I_G$. The latter property is also known to be equivalent to the following  {\sl spectral gap
condition\/}:  there exist a positive lower bound for  
the spectrum of the Laplacian $\bold
\Delta$ on $K\backslash \ggm$, where $K$ is a maximal compact subgroup
of $G$. 
 

If $G$ is a direct product of simple
groups $G_1,\dots,G_l$, one can decompose $\bold
\Delta$ as a sum $\bold
\Delta_1 + \dots + \bold \Delta_l$, where $\bold \Delta_i$ corresponds
to coordinates coming from $G_i$. Then a lower bound for
the spectrum of  $\bold \Delta_i$ amounts to the 
 isolation of $\rho_0|_{G_i}$  from
the trivial representation $I|_{G_i}$ of $G_i$. In the paper \cite{KM} it was
implicitly conjectured that restrictions $\rho_0|_{G_i}$ are isolated from
$I|_{G_i}$. Theorem 1.12 proves this conjecture in the non-uniform
lattice case. The main ingredient of the proof is an explicit bound
for the bottom of spectra of Laplacians given by M.-F.~Vigneras in
\cite{V}. The reduction to the case where these bounds are applicable
is based on Lemma 3.1, the Arithmeticity Theorem and the restriction
technique of 
M.~Burger and P.~Sarnak. We now present the 

\demo{Proof of Theorem 1.12} If $G$ is
simple, the claim follows from \cite{B, Lemma 4.1}. Therefore we can
assume  that the $\br$-rank of $G$ is greater than $1$. By Margulis'
Arithmeticity Theorem (see \cite{Z, Theorem 6.1.2} or \cite{M, Chapter
IX}), $\Gamma$ is an 
arithmetic lattice in $G$. That 
is, there exists a semisimple algebraic $\bq$-group $\bold G$ and a
surjective  homomorphism $\ph:\bold G(\br)^0\mapsto G$ such that:


\roster
\item"{(i)}" $\text{Ker}\,\ph$ is compact, and

\item"{(ii)}" the subgroups $\ph\big(\bold G(\bz) \cap \bold G(\br)^0\big)$ and $\Gamma$ are commensurable. 
\endroster

Further, since $\Gamma$ is
non-uniform and $G$ is center-free, $\bold G$ can be taken to be
connected and adjoint, and
$\text{Ker}\,\ph$ to be trivial (see \cite{Z, Corollary 6.1.10}). By
(ii) above, the spaces $\ggm$ and $G/\ph\big(\bold G(\bz)\big)$ have a
common finite covering. In view of Lemma  3.1, without loss of
generality one can 
assume that $\Gamma = \bold G(\bz)$ and $G = \bold G(\br)$. 

Write $\bold G$ in the form $R_{k/\bq}\tilde \bold G$, where $k$ is a
finite extension of $\bq$, $\tilde \bold G$ is an absolutely
$k$-simple $k$-group, and $R_{k/\bq}$ stands for Weil's restriction of
scalars functor (see \cite{T2, 3.1.2}). Namely, $\bold G = \prod_{i =
1}^l\tilde \bold 
G^{\sigma_i}$, where $\sigma_1,\dots,\sigma_l$ are 
distinct  imbeddings of $k$ into $\bc$. 
This way, factors $G_i$ of
$G$ can be identified with $\tilde \bold
G^{\sigma_i}(\br)$ if $\sigma_i$ is real, or with $\tilde \bold
G^{\sigma_i}(\bc)$ if $\sigma_i$ is complex.

 Since $\Gamma$ is non-uniform, $\tilde\bold G$ is isotropic over $k$ (see \cite{M, Theorem I.3.2.4(b)}),
therefore (see \cite{T1, 3.1, Proposition 13} or \cite{M, Proposition
I.1.6.3}) there exists  
a $k$-morphism $\tilde \alpha:\bold{S}\bold{L}_2\mapsto  \tilde\bold G$
with finite 
kernel. 
Denote the $\tilde \alpha$-image of $\bold{S}\bold{L}_2$ by $\tilde
\bold H$, and let $\bold H = R_{k/\bq}\tilde \bold H$
and $H = \bold H(\br)$.  Clearly to show that $\rho_0|_{G_i}$ is
isolated from $I|_{G_i}$, it will be enough to prove that  $\rho_0|_{H_i}$ is
isolated from $I|_{H_i}$, where $H_i$ are almost simple factors of $H$,
isomorphic to $\tilde \bold H^{\sigma_i}(\bc)$ for complex imbeddings
$\sigma_i$ and to $\tilde \bold H^{\sigma_i}(\br)$ for real ones.

We now use Theorem 1.1 from the paper \cite{BS}, which guarantees that
$\rho_0|_{H}$ lies in the closure of the {\sl automorphic
spectrum\/} of $H$ (the latter stands for irreducible
components of representations of
$H$ on all the spaces $L^2(H/\Lambda)$ where
$\Lambda$ is a congruence subgroup of $\bold H(\bz)$). Denote by
$\bold L$ the algebraic group $R_{k/\bq}\bold{S}\bold{L}_2$
and by 
$\alpha$ the isogeny $\bold{L} \mapsto \bold H$
induced by $\tilde \alpha$. Note that \hs s
$H/\Lambda$ can be identified with
$\bold{L}(\br)/\alpha^{-1}(\Lambda)$, and preimages
of congruence subgroups of $\bold H(\bz)$ are congruence subgroups of
$\bold{L}(\bz)$. 
Therefore it suffices to check that nontrivial irreducible components
of regular representations of almost $\br$-simple factors of
$\bold{L}(\br)$ 
on   $L^2\big(\bold{L}(\br)/\Lambda\big)$ are
uniformly isolated from the trivial representation for all $i =
1,\dots,r$ and all  principal congruence 
subgroups $\Lambda$ of $\bold{L}(\bz)$. The latter
statement is a 
reformulation 
 of one of the corollaries in Section VI of the paper
\cite{V}, with the uniform bound for the first nonzero eigenvalue of
the corresponding Laplace operators being equal to $3/16$ for real and
$3/4$ for complex imbeddings $\sigma_i$. 
\qed\enddemo


\example{3.3. Remark} One can also  prove Theorem 1.12 without using
Lemma  3.1 by extending the result of Vigneras to arbitrary subgroups
of  $\bold H(\bz)$ rather than congruence subgroups.  For this one can
use the centrality of the congruence kernel for higher rank groups, see \cite{R2},
and the results of Y.~Flicker \cite{F} on lifting of automorphic representations to
metaplectic coverings of $GL_2$. This way it should be possible to get
an explicit uniform (in all $G$ and $\Gamma$) bound for the
neighborhood of the trivial representation which is disjoint from all
the restrictions $\rho_0|_{G_i}$.
\endexample

\subhead{3.4}\endsubhead We now turn to the paper \cite{KM}, where 
the well-known (from the work of Harish Chandra, Howe, Cowling and
Katok-Spatzier) connection between isolation properties of $\rho_0$ and
exponential decay of its matrix coefficients has been made explicit. Let
$G$ be a connected semisimple Lie group with finite center, $K$ its
maximal compact subgroup.  
Take an orthonormal basis $\{X_i\}$ of the Lie algebra of $K$, and
denote by $\w$  the differential operator $1-\sum_{i =
1}^{\text{dim}(K)} X_i^2$ (see \cite{W, \S 4.4.2}). 

\proclaim{Theorem \rm (see \cite{KM, Corollary 2.4.4} and a
correction on p.~172)} Let $\Pi$ be a family of 
unitary representations of $G$ such that the restriction of $\Pi$ to any simple factor of $G$ is isolated
 from the trivial representation. 
Then there
exist a universal constant $B>0$, a positive integer   $l$
(dependent only on $G$) and $\beta>0$ (dependent on $\Pi$ and on the
choice of the bi-$K$-invariant norm $\|g\| = \text{\rm dist}(g,e)$  on $G$) 
such that for any $\rho\in\Pi$, any $C^\infty$-vectors $v$, $w$ in a
representation space of $\rho$   
and any $g\in G$  one has
$$
\left|\big(\rho(g)v,w\big)\right|\le
Be^{-\beta\|g\|}\|\w ^l(v)\|\|\w ^l(w)\|\,.\tag 3.2
$$
\endproclaim

Combining
Theorem 3.4 and Theorem 1.12, we obtain  
the following

\proclaim{3.5. Corollary}  
 Let $G$ be a connected semisimple center-free Lie
group  without compact factors,  $\Gamma$ an irreducible
non-uniform lattice in $G$, $X = \ggm$, $\mu$ the normalized Haar
measure on $X$. 
Then there
exist  $B,\beta>0$ and   $l\in \bn$   
such that  for any two functions $\ph, \psi \in C^\infty(X) \cap L^2(X)$    
and any $g\in G$  one has
$$
\big|(g\ph,\psi) - \mu(\ph) \mu(\psi)\big|\le
Be^{-\beta\|{g}\|}\|\w ^l(\ph)\|\|\w ^l(\psi)\|\,.\tag 3.2
$$
\endproclaim

\demo{Proof} The family $\Pi = \{\rho_0\}$ satisfies the assumption of
Theorem 3.4 in view of Theorem 1.12. Therefore one can apply (3.2) to
the functions $\ph - \mu(\ph)$ and $\psi - \mu(\psi)$. \qed\enddemo

 \heading{\S 4. A quantitative version of Theorem 1.8}\endheading

\subhead{4.1}\endsubhead 
Let  $G$, $\Gamma$ and
$\mu$ be as in Theorem 1.12,  and denote the (noncompact) \hs\ $\ggm$ by $X$. 
Our first goal is to apply Lemma
2.6 to certain sequences of functions on $X$. For $l\in\bn$ and $C >
0$, say that 
$h \in C^\infty(X) \cap L^2(X)$ is {\sl $(C,l)$-regular\/} if
$$
\|\w ^l(h)\| \le C\cdot \mu(h)\,.
$$


\proclaim{Proposition} Assume that 
$F=
\{f_{t}\mid {t}\in \bn\}$ is an ED sequence of elements of $G$. Take
$l\in\bn$ as in Corollary 3.5 and an arbitrary $C > 0$, and let  $\ch
= \{h_t\}$
be a  
sequence of   $(C,l)$-regular functions on $X$ such that {\rm (2.7)} holds. Then the
twisted sequence $\ch^F$
satisfies {\rm (SP)}; in particular, {\rm (2.8)} holds and 
$$
\lim_{N\to\infty}\dfrac{S_{\ch^F,N}(x)}{E_{\ch,N}} = 1 \quad\text{for
$\mu$-a.e.\ 
$x\in X$}
$$
whenever $\ch$ is nonsummable.
 \endproclaim 

\demo{Proof} In view of (2.9), one has to estimate the sum
$$
 \sum_{s,t = M}^{N} \big((f_s^{-1}h_s,f_t^{-1}h_t) -
\mu(h_s)\mu(h_t)\big)\tag 4.1
$$
 from above. Observe that, since $\mu$ is $F$-invariant,
$(f_s^{-1}h_s,f_t^{-1}h_t) -
\mu(h_s)\mu(h_t)$ is equal to 
$$
\split
 (h_s,f_sf_t^{-1}h_t) -
\mu(h_s)\mu(h_t)  \un{(by Corollary 3.5)}{\le}
  &Be^{-\beta\|{f_sf_t^{-1}}\|} \|\w ^l(h_s)\|\|\w ^l(h_t)\| \\
\un{(by
the 
$(C,l)$-regularity of $h_s$, $h_t$)}{\le}
BC^2e^{-\beta\|{f_sf_t^{-1}}\|}\mu(h_s)\mu(h_t) \un{(by (2.7))}\le
&BC^2e^{-\beta\|{f_sf_t^{-1}}\|}\mu(h_t)\,. 
\endsplit
$$
Therefore the sum (4.1) is not bigger than 
$$
BC^2 \sum_{s,t = M}^N e^{-\beta\|{f_sf_t^{-1}}\|}\mu(h_t) = BC^2
\sum_{t =M}^N \mu(h_t) \sum_{s =M}^N e^{-\beta\|{f_sf_t^{-1}}\|}
\le BC^2\cdot\sup_{t\in\bn}\sum_{s=1}^\infty
e^{-\beta\|f_sf_t^{-1}\|} \cdot E_N 
\,.
$$
In view of (ED), the constant in the right hand side is finite, and
(SP) follows; the ``in particular'' part is then immediate from Lemma 2.6.
\qed \enddemo

\subhead{4.2}\endsubhead  Let  now $\Delta$ be a
DL function on $X$. Similarly to (1.4), for $z \in\br$ we will denote
by $A(z)$ the set $\{x\in X\mid 
\Delta(x)\ge z\}$ (note that it follows from (DL) that $A(z)$ is never
empty). To prove a quantitative strengthening of Theorem 
1.8 that we are after, we need to learn how to
approximate the sets $A(z)$  by smooth functions. 
 
\proclaim{Lemma} Let  $\Delta$ be a
DL function on $X$. Then for any $l\in\bn$ there exists  $C > 0$ such
that for every $z \in\br$ one can find two  $(C,l)$-regular nonnegative
functions $h'$ and $h''$ on $X$ such that  
$$
h' \le 1_{A(z)} \le h''\quad \text{and}\quad c\cdot\mu\big(A(z)\big) \le 
\mu(h') \le \mu(h'') \le \frac1c\mu\big(A(z)\big)\,,\tag 4.2
$$
with $c$ as in {\rm (DL)}.
  \endproclaim

\demo{Proof}  For  $\vre  > 0$, let us
denote by $A'(z,\vre )$ the set of all points of $A(z)$ which are not
$\vre $-close to $\partial A(z)$, i.e.~$A'(z,\vre )\df \{x\in A(z)\mid
\text{dist}\big(x,\partial A(z)\big) \ge \vre \}$, and  by $A''(z,\vre )$ the
$\vre $-neighborhood of  $A(z)$, i.e.~$A''(z,\vre )\df \{x\in X \mid \text{dist}\big(x,
A(z)\big) \le \vre \}$.  (If $A(z) = X$, the above sets will 
coincide with $X$.) 

Choose $\delta$ and $c$ according to (DL).   Then, using the
uniform continuity of $\Delta$, find $\vre  > 0$ such that 
$$
|\Delta(x) -
\Delta(y)| < \delta\text{ whenever dist}(x,y) < \vre \,.\tag 4.3
$$
 From (4.3) it immediately follows that for all $z$ one has 
$
A(z+\delta) \subset A'(z,\vre )\subset A''(z,\vre )\subset
A(z-\delta)$, therefore  one
can apply (DL) to conclude that  
$$
c\cdot\mu\big( A(z)\big)\le \mu\big( A'(z,\vre )\big) \le \mu\big(
A''(z,\vre )\big) \le \frac1c\mu\big(A(z)\big)\,.\tag 4.4 
$$

Now  take a nonnegative $\psi\in C^\infty(G)$ of $L^1$-norm 1 such that
supp$(\psi)$ belongs to the ball of radius $\vre /4$ centered in $e\in G$.
 Fix  $z\in\br$ and consider   functions $h' \df \psi* 1_{A'(z,\vre /2)}$ and $h'' \df \psi* 1_{A''(z,\vre /2)}$. Then one clearly has
$$
1_{A'(z,\vre )} \le h' \le 1_{A(z)} \le h'' \le 1_{A''(z,\vre )} \,,
$$
which, together with  (4.4),   
 immediately implies (4.2). It remains to observe that  
$\|\w^l h'\| =  \|\w^l(\psi* 1_{A'(z,\vre /2)})\|  = \|\w^l(\psi)*
1_{A'(z,\vre /2)}\|$, so by the Young
inequality, 
$$
\|\w^l h'\| \le \|\w^l(\psi)\| \cdot \mu\big(A'(z,\vre /2)\big)  \le
\|\w^l(\psi)\|  \cdot \mu\big(A(z)\big)   \un{(4.2)}\le
\frac1{c}\|\w^l(\psi)\| \mu(h')\quad\text{for any }l\in
\bn\,.
$$  
 Similarly $\|\w^l h''\| \le \|\w^l(\psi)\| \cdot
\mu\big(A''(z,\vre /2)\big) \un{(4.4)}\le \|\w^l(\psi)\| \cdot
\frac1{c}\mu\big(A(z)\big) \le \frac1{c}\|\w^l(\psi)\|  \cdot\mu(h'')$,
hence, with $C = \frac1{c}\|\w^l(\psi)\|$, both  $h'$ and  $h''$ are
$(C,l)$-regular, and the lemma is proven. \qed\enddemo 


\subhead{4.3}\endsubhead We now state and prove the promised quantitative
strengthening of Theorem 1.8.

\proclaim{Theorem}  Let $G$,  $\Gamma$, $F =
\{f_{t}\}$ and  $\Delta$ be as in Theorem 1.8, 
and let $\{r_{t}\}$ be
a sequence of real numbers such that  
$$
{\dsize\sum_{{t}=1}^{\infty} \td(r_{t})} = \infty\,.\tag 4.5
$$
Then for some positive 
$c \le 1$ and for  almost
all $x\in \ggm$ one has
$$
 c \le \liminf_{N\to\infty}\frac {\#\{ 1 \le {t}\le N\mid \Delta(f_{t}x)
\ge r_{t}\}}  {\sum_{t=1}^{N} \td\big(r_{t}\big)} \le \limsup_{N\to\infty}\frac {\#\{ 1 \le {t}\le N\mid \Delta(f_{t}x)
\ge r_{t}\}}  {\sum_{t=1}^{N} \td\big(r_{t}\big)} \le \frac1c\,.
$$
\endproclaim

It is clear that Theorem 1.8 is a direct consequence of Lemma 2.3 and
the first of 
the above inequalities. Note that D.~Sullivan proved that in the
setting of Theorem 1.2 one has a positive lower bound for 
$$ 
\limsup_{N\to\infty}\frac{\#\{ 1 \le
{t}\le N\mid  
\text{\rm dist}\big(y_0,\gamma_{t}(y,\xi)\big) \ge  r_{t}\}}
{\sum_{t=1}^{N} e^{-{k} r_{t}}}
$$  
for  almost
all  $\xi\in S_y(Y)$  (see \cite{Su, \S 9, Remark (2)}). 

\demo{Proof}  First let us assume that the center of $G$ is trivial;
after that we will reduce the general case to the center-free
situation. Observe that from the existence of a DL function
$\Delta$ 
on $X$ it follows that $X$ can not be compact: indeed, $\Delta$ must
be uniformly continuous, but unbounded in
view of (DL). Take $l$ as in Corollary 3.5 and $C$ as in Lemma 4.2. 
For any ${t}\in\bn$, let $h_{t}'$ and $h_{t}''$ stand for the
$(C,l)$-regular functions which one associates with  the set
$A(r_{t}) = \{x\in X\mid 
\Delta(x)\ge r_{t}\}$  by means of
Lemma 4.2, and let us denote 
$$
\ca = \{A(r_{t})\mid {t}\in \bn\}\,,\quad \ch' = \{h'_{t}\mid {t}\in
\bn\}\,,\quad \ch'' = \{h''_{t}\mid {t}\in \bn\}\,.
 $$
By (4.5), the sequence $\ca$ is 
nonsummable; hence, in view of (4.2), the same can be said 
about $\ch'$ and $\ch''$.  Also it is clear from the construction that
$ \mu(h'_{t}) 
\le \mu(h''_{t}) \le 1$. Therefore,  as $N\to\infty$, by Proposition
4.1  the ratios  
${S_{(\ch')^F, N}(x)}/{E_{\ch', N}}$ and ${S_{(\ch'')^F,
N}(x)}/{E_{\ch'', N}}$ tend to $1$ for {$\mu$-a.e.}~$x\in X$.
But from (4.2) it follows that $ S_{(\ch')^F, N} \le S_{\ca^F, N} \le
S_{(\ch'')^F, N}$ and $\frac1c{E_{\ch', N}}\ge E_{\ca, N} \ge c
\cdot {E_{\ch'', N}}$ for any $N\in\bn$. Therefore $\mu$-almost
everywhere one has  
$$
 c =  \lim_{N\to\infty}\dfrac {S_{(\ch')^F, N}(x)}{\frac1c E_{\ch',
N}} \le \liminf_{N\to\infty}\dfrac  {S_{\ca^F, N}(x)}{E_{\ca, N}} \le
\limsup_{N\to\infty}\dfrac  {S_{\ca^F, N}(x)}{E_{\ca, N}} \le
\lim_{N\to\infty}\dfrac {S_{(\ch'')^F, N}(x)}{c \cdot E_{\ch'', 
N}} = \frac1c\,,
$$
and the statement of the theorem follows.

Now let us look at what happens if $G$ has nontrivial center $Z$. Let us denote the
quotient group $G/Z$ by $G'$, the homomorphism
$G\mapsto G'$ by $p$,  and the induced map
$X\mapsto X'\df G'/p(\Gamma)$ by $\bar p$.
Since $\Gamma  Z$ is discrete \cite{R1, Corollary 5.17},
$p(\Gamma)$ is also discrete, hence $Z/(\Gamma\cap Z)$ is finite. This means that 
$(X,\bar p)$ is a finite covering of $X'$; moreover, one can choose
representatives $g_1,\dots,g_l$ ($g_1 = e$) from cosets of $Z/(\Gamma\cap Z)$
which will act isometrically on $X$. In particular, the distance between
$x\in X$ and $g_ix$, $1\le i \le l$, is uniformly bounded by some
constant $C$. Now,  given a DL
function $\Delta$ on $X$, define $\Delta'$ on $X'$ by $\Delta'\big(\bar p(x)\big)
= \frac1l\sum_{y\in \bar p^{-1}\circ \bar p (x)} \Delta(y)=\frac1l \sum_{i = 1}^l
\Delta(g_ix)$. Then from the uniform continuity of $\Delta$ it follows
that $\Delta'$ is also uniformly  continuous, and for some constant
$C'$ one has
$$
|\Delta'\big(\bar
p(x)\big) -  \Delta(x)| \le C'\quad\forall\,x\in X\,.\tag 4.6
 $$
 Therefore for any $z > 0$,
$\Phi_{\Delta'}(z) = \mu\big(\big\{x\in X\mid \Delta'\big(\bar
p(x)\big)\ge z\big\}\big)$ is bounded between  $\Phi_{\Delta}(z+C')$
and  $\Phi_{\Delta}(z-C')$. This implies that $\Delta'$ satisfies
(DL) as well; moreover, 
$$
\frac{\Phi_{\Delta'}(z)}
{\Phi_{\Delta}(z)}\text{ is uniformly bounded between two positive
constants. }\tag 4.7
$$

Finally,  assume that (4.5) holds  and $F\subset G$
is ED. It follows that $\{p(F)\}$ is also ED, and from
(4.7) one 
deduces that $
{\sum_{{t}=1}^{\infty} \Phi_{\Delta'}(r_{t})} =
\infty$ as well.  Therefore one can  use the center-free case of
Theorem 4.3 and $\Delta'$ being a DL function to conclude that  for some
$0 < c \le 1$ and for  $\mu$-almost
all $x\in X$ one has
$$
 c \le \liminf_{N\to\infty}\dfrac {\#\{ 1 \le {t}\le N\mid \Delta'\big(\bar p(f_{t}x)\big)
\ge r_{t} + C'\}}  {\sum_{{t}=1}^{N} \Phi_{\Delta'}(r_{t})}
$$
and
$$
 \limsup_{N\to\infty}\dfrac {\#\{ 1 \le {t}\le N\mid \Delta'\big(\bar p(f_{t}x)\big)
\ge r_{t} - C'\}}  {\sum_{{t}=1}^{N} \Phi_{\Delta' }(r_{t})} \le \frac1c\,.
$$
Clearly  (4.6) implies that 
$$
\Delta'\big(\bar p(f_{t}x)\big)
\ge r_{t} + C' \thus \Delta(f_{t}x)
\ge r_{t} \thus \Delta'\big(\bar p(f_{t}x)\big)
\ge r_{t} - C'\,.
$$
 Therefore to finish the proof it remains to replace
the values of $\Phi_{\Delta'}$ by those of $\td$, sacrificing no more
than a multiplicative constant in view of (4.7).
 \qed\enddemo

\demo{{\bf 4.4.}\ \ Proof of
Theorems 1.7 and 1.9} Recall that in part (a) of Theorem 1.9 we are
given  a sequence 
$F = \{f_{t}\}  
= \{\exp({\vz_{t}})\}$ such that (1.6) holds. It is easy  to check that
$F$ 
satisfies (ED):  for any $\beta >
0$ one can write
$$
\split
\sup_{t\in\bn}\sum_{s=1}^\infty e^{-\beta\|f_sf_t^{-1}\|} &=
\sup_{t\in\bn}\sum_{s=1}^\infty e^{-\beta\|\vz_s-\vz_t\|} 
\le 
\sup_{t\in\bn}\sum_{n=0}^\infty e^{-\beta n}\#\{s\mid n \le
\|\vz_s-\vz_t\| \le n + 1\} \\ &\un{(1.6)}{\le} \text{const}\cdot \sum_{n=0}^\infty
n^{\text{dim}(\a)}e^{-\beta n}  < \infty\,.
\endsplit
$$
Therefore Theorem 
1.8 applies and one concludes that $\cb(\Delta)$ is Borel-Cantelli for
$F$. Part 
(b) is then immediate from 
Corollary 2.4 and Lemma 2.8. 
It remains to notice that  Theorem 1.7 is a
special case of Theorem 1.9, with  $\vz_{t} = {t}\vz$, $d =
1$, $\d =
\br \vz$ and $\d_{\sssize +} = \{t \vz\mid t\ge 0\}$. 
\qed\enddemo


 \heading{\S 5. Distance functions are DL}\endheading

\subhead{5.1}\endsubhead The goal of the section is to prove the
following 

\proclaim{Proposition}   Let $G$ be a connected semisimple 
Lie
group, 
$\Gamma$ a non-uniform irreducible\footnote{Again, the
proposition  is also true for reducible lattices, see
\S 10.2.} 
lattice in $G$, $K$ a maximal compact subgroup of $G$, $\mu$ the
normalized Haar measure on $\ggm$, $x_0$ a point in 
$\ggm$,   {\rm
dist}$(\cdot,\cdot)$
 a Riemannian metric  on  $\ggm$ 
chosen by fixing  a 
right invariant 
Riemannian metric 
on $G$ bi-invariant with respect to
$K$. Then there exists $k > 0$ such that 
the function  $\text{\rm
dist}(x_0,\cdot)$ is   $k$-DL. 
 \endproclaim

\example{5.2. Remark}  Let  $(X_1,x_1)$  and $(X_2,x_2)$ be pointed
metric spaces with probability measures $\mu_1$ and $\mu_2$,  and let  $\pi:X_1\mapsto X_2$ be a
measurable surjective  map which almost preserves distances from base
points (i.e.~with $\sup_{x\in X_1}\big|\text{\rm
dist}(x_1,x) - \text{\rm
dist}\big(x_2,\pi(x)\big)\big| < \infty$) and satisfies the following
property: for
some positive $c<1$ one has 
$$
c \mu_2(A) \le \mu_1\big(\pi^{-1}(A)\big) \le \frac1c
\mu_2(A)\quad\text{for any }A\subset X_2\,.
$$
Then the function  $\text{\rm
dist}(x_1,\cdot)$ on $X_1$ is   $k$-DL iff so is $\text{\rm
dist}(x_2,\cdot)$ on $X_2$. 
 This observation will
be
used many times in the proof below, sometimes without explicit
mention. Examples include:

\roster 
\item"{$\bullet$}" $X_1 = X_2$, $\mu_1 = \mu_2$ (shift of base point);

\item"{$\bullet$}" $X_1 = X_2 \times Q$ (the direct product of metric
and probability spaces), $\pi$ the projection on $X_2$, diam$(Q)<
\infty$;

\item"{$\bullet$}" $X_1\overset{\pi}\to\mapsto 
X_2$ a finite covering,
$\mu_2 = \pi(\mu_1)$.
 \endroster
\endexample

\demo{{\bf 5.3.}\ \ 
Proof of Proposition 5.1}  First suppose that the 
$\br$-rank of $G$ is greater than 
$1$. Then, using the Arithmeticity Theorem, as in the
proof of Theorem 
1.12 (see \S 3.2) we can assume that  
$G =
\bold G(\br)$, where  $\bold G$
is a semisimple algebraic $\bq$-group and $\Gamma$ is an arithmetic
subgroup of $\bold G(\bq)$. 

We now need to use the reduction theory for arithmetic groups. Let 
$\bold T$ be a maximal 
$\bq$-split
torus of $\bold G$. 
Denote the identity
component of $\bold 
T(\br)$ by $A$, and its Lie algebra by $\a$. Let  $\Phi$ be a system of
$\bq$-roots associated with 
$\a$. Choose an ordering of  $\Phi$, let  $\Phi^{\sssize +}$ (resp.~$\Phi^s$)
be the set of positive (resp.~simple) roots, and let ${\a_{\sssize +}}$
stand for the closed $\bq$-Weyl chamber in $\a$ defined by 
$
{\a_{\sssize +}} \df\{\vz\in \a\mid \alpha(\vz )\ge 0\ \forall\,\alpha
\in\Phi^s\}\,.$
We set ${A_{\sssize +}} \df \exp({\a_{\sssize +}})$. 

Let $G = KAMU$ be a (generalized) Iwasawa decomposition for $G$,
where $K$ is maximal compact in $G$, $U$ is unipotent and $M$ is
reductive (here $A$ centralizes $M$ and normalizes $U$). Then one
defines a {\sl generalized Siegel set\/} $\Cal S_{Q,\tau}$ as follows:
$\Cal S_{Q,\tau}\df K\exp(\a_\tau) Q$, where $Q$ is relatively compact in
$MU$, $\tau \in\br$ and ${\a_\tau} \df\{\vz\in \a\mid \alpha(\vz)\ge \tau\ \forall\,\alpha
\in\Phi^s\}$. 
%
It is known that a
finite union of translates of such a set (for suitable $Q$ and $\tau$)
forms a {\sl weak fundamental set\/} for the $G$-action on
$\ggm $.  More
precisely, the following is true:  

\proclaim{5.4. Theorem \rm(\cite{Bo, \S13} or \cite{L, Proposition
2.2})}  Let  $\bold G$
be a semisimple algebraic $\bq$-group and $\Gamma$  an arithmetic
subgroup of $\bold G(\bq)$. Then  there exist a  generalized Siegel
set $\Cal S = 
\Cal S_{Q,\tau}\subset G = \bold G(\br)$ and $\{q_1,\dots,q_m\}\subset \bold G(\bq)$ such that
the  union $\Omega \df \cup_{i = 1}^m 
\Cal Sq_i$ satisfies the following two properties:
\roster
\item"{(i)}" $G = \Omega \Gamma $;

\item"{(ii)}" for any $q\in \bold G(\bq)$, the set $\{\gamma\in
\Gamma\mid \Omega q\cap \Omega \gamma\}$ is finite.
\endroster
 \endproclaim

In other words, the restriction to $\Omega$  of the natural projection
$\pi$ of $G$   
onto $\ggm $ is surjective and at most finite-to-one.

We now want to study metric properties of the restriction
$\pi|_\Omega$. Since the distance on $\ggm $ is defined by
dist$_{\ggm }\big(\pi(x),\pi(y)\big) =
\inf_{\gamma\in\Gamma}\text{dist}_G(x,y\gamma)$, one clearly has
dist$_{\ggm }\big(\pi(x),\pi(y)\big)$ $\le \text{dist}_G(x,y)$ for
any $x,y \in G$. The converse estimate, with $x,y$ taken from a
Siegel set, has been known as Siegel's Conjecture. 
Its proof is due to
J.~Ding for $G = SL_n(\br)$ and to E.~Leuzinger and L.~Ji
 (independently) for the general case. Specifically, the following
statement has been proved:

\proclaim{5.5. Theorem \rm(\cite{L, Theorem 5.7} or \cite{J, Theorem
7.6})} For  $\bold G$, $\Gamma$, $\Cal S$ and $\{q_1,\dots,q_m\}$ as
in Theorem 5.4, there exists a positive constant $D$ such that 
$$
\text{\rm dist}_{G}\big(xq_i,yq_j\gamma\big) \ge \text{\rm
dist}_G(x,y) - D
$$
for all $i,j = 1,\dots,m$, $\gamma\in\Gamma$ and $x,y\in \Cal S$.
 \endproclaim

In view of the last two theorems and Remark 5.2, it is enough to prove
that the function $\text{\rm
dist}_G(x_0,\cdot)$ on $\Omega$ is  $k$-DL for some $k > 0$
and $x_0\in\Omega$ (with respect to suitably scaled Haar
measure). Further, since  the metric on $G$ is right 
invariant, it suffices to consider just one copy $\Cal S =
K\exp(\a_\tau) Q$ of the 
Siegel set instead of the union $\Omega$ of several translates
thereof. 

Our next goal is to reduce the problem to the restriction of the
distance function to $\exp(\a_\tau)$. Since  the metric on $G$ is right 
invariant and bi-$K$-invariant, the projection $G = KAMU \mapsto A$ is almost
distance preserving (in the sense of Remark 5.2). Furthermore, cf.~\cite{Bou1,
Ch.~VII, \S 9, Proposition 13}, the Haar
measure on $G$ is being sent to the measure $\delta(a)\,da$, where $da$
is a Haar measure on $A$ and $\delta$ is the restriction of the modular
function of the group $AMU$ to $A$. Put differently, $\delta(a)$ is the
modulus of the automorphism $x\mapsto axa^{-1}$ of $MU$ (equivalently, of
$U$, since $M$ is centralized by $A$). Therefore, if $a = \exp(\vz )$,
$\vz \in\a$,  the value of $\delta$ at $a$ is equal to $e^{\text{tr}(-\text{ad}\,\vz )} =
e^{-\rho(\vz )}$, where $\rho \df \sum_{\alpha\in\Phi^{\sssize +}}\alpha$
is the sum of the positive roots. Since the metric on $A$ is
carried from $\a$ by the exponential map, it suffices to find $k$ such
that 
the function $\vz\mapsto\|\vz\|$ on ${\a_\tau}$ (equivalently, on
${\a_{\sssize +}}$, since ${\a_\tau}$  is an isometric
translate of the latter) is  $k$-DL
with respect to the measure const$\cdot
e^{-\rho(\vz )}\,d\vz $. 

Let $\{\alpha_1,\dots,\alpha_n\}$ be the simple roots, and
$\{\omega_1,\dots,\omega_n\}$ the dual system of fundamental
weights (that is, with $\alpha_i(\omega_j) = \delta_{ij}$, $i,j =
1,\dots,n$). One can write 
$$
\rho = \sum_{i = 1}^n k_i\alpha_i\,,\tag 5.1
$$
where $k_i$ are positive integers. The following lemma is what one
needs to complete the proof:

\proclaim{5.6. Lemma} Let $k = \min_{i = 1,\dots,n}
\frac{k_i}{\|\omega_i\|}$. Then there exist $C_1,C_2 > 0$ such that
for any $z > 0$, the ratio of  $\int_{\{\vz\in{\a_{\sssize +}},\, \|\vz\| \ge z\}}e^{-\rho(\vz )}\,d\vz $ and $e^{-kz}$ is bounded between $C_1$ and $C_2$.
 \endproclaim

\demo{Proof} Without loss of generality assume that
$\frac{k_i}{\|\omega_i\|}$ is equal to $k$ for $1\le i \le j$ and is
strictly bigger than $k$ for $i> j$. For $r > 0$, denote by $\Sigma_r$ the
intersection of ${\a_{\sssize +}}$ and the 
sphere of radius $r$ centered at the origin. This is a spherical
simplex with extremal points (vertices) given by $\vz _i\df\frac
r{\|\omega_i\|}\omega_i$. From the strict convexity of the ball it follows
that $\rho|_{\Sigma_r}$ attains its minimal value $kr$ at the points $\vz_i$, $1\le i \le j$. Furthermore, one can choose 
$\vre,\vre',c > 0$ such that uniformly in $r > 0$ the set
$$
\Sigma_{r,\vre} \df \{\vz\in \Sigma_r ,\, \rho(\vz) \le (k + \vre) r\}
$$
belongs to the union of $\vre'r$-neighborhoods of the points  $\vz_i$,
$1\le i \le j$, and on each of these neighborhoods one has $\rho(\vz)
- kr \ge c\|\vz - \vz_i\|$. 

Denote by
$\sigma$ the induced Lebesgue measure on  $\Sigma_r$. Clearly to
establish the desired 
upper estimate for 
$$
\int_{\{\vz\in{\a_{\sssize +}} ,\, \|\vz\| \ge
z\}}e^{-\rho(\vz )}\,d\vz =  \int_z^\infty \int_{\Sigma_r}
e^{-\rho(\vz )}\,d\sigma(\vz)\,dr
$$ 
it suffices to prove
that
$\int_{\Sigma_r} e^{-\rho(\vz )}\,d\sigma(\vz)$ is not greater
than\footnote{The values of constants  
in the proof below are
independent on $r$.}
const$\cdot e^{-kr}$. The latter inequality follows since 
$$
\split 
\int_{\Sigma_r} e^{-\rho(\vz )}\,d\sigma(\vz)  &\le \int_{\Sigma_r\ssm
\Sigma_{r,\vre}} e^{-\rho(\vz )}
\,d\sigma(\vz)  + \int_{\Sigma_{r,\vre}} e^{-\rho(\vz )}
\,d\sigma(\vz)\\
 &\le \int_{\Sigma_r} 
e^{-(k + \vre) r}\,d\sigma(\vz)  + \sum_{i = 1}^j\int_{\{\vz\in
\Sigma_r ,\, \|\vz - 
\vz_i\| \le \vre' r\}} e^{-(kr + c\|\vz - \vz_i\|)}\,d\sigma(\vz)\\&\le
\text{const}\cdot r^{n-1}e^{-(k + \vre) r}  + 
\text{const}\cdot e^{-kr}\int_{\br^{n-1} } e^{-c\|\vx\|}\,d\vx \le
\text{const}\cdot e^{-kr}\,.
\endsplit
$$

As for the lower estimate, the set $\{\vz\in{\a_{\sssize +}} ,\, \|\vz\| \ge
z\}$ clearly 
contains
the translate $\vz_1 + {\a_{\sssize +}}$ of 
${\a_{\sssize +}}$, where, as before, $\vz_1 =  \frac
z{\|\omega_1\|}\omega_1$ and $\rho(\vz_1) = kz$; therefore  
$$ 
\int_{\{\vz\in{\a_{\sssize +}} ,\, \|\vz\| \ge
z\}} e^{-\rho(\vz )}\,d\vz \ge \int_{\vz_1 + {\a_{\sssize +}}}
e^{-\rho(\vz )}\,d\vz 
= \int_{{\a_{\sssize +}}} e^{-\rho(\vz + \vz_1 )}\,d\vz =
e^{-kz}\int_{\a_{\sssize +}} 
e^{-\rho(\vz )}\,d\vz \,,
$$
which finishes the proof. \qed\enddemo 

To complete the proof of Proposition 5.1 it remains to observe that in
the case when the $\br$-rank of $G$ is equal to 
$1$, the proof can be written along the same lines, by means of  the
description \cite{GR} of fundamental domains for lattices in rank-one groups. 
 \qed\enddemo

\subhead{5.7}\endsubhead Note that the above proof, via Lemma 5.6,
provides a constructive way 
to express the exponent $k$  for any homogeneous space $\ggm$ via
parameters of the corresponding system $\Phi$ of $\bq$-roots. For
example, if $G = SL_n(\br)$ and the metric on $G$
is given by the Killing form, one can compute (using e.g.~\cite{Bou2,
Planche I}) the norms of fundamental weights
$\omega_1,\dots,\omega_{n-1}$: 
$$
\|\omega_i\|^2 = \frac{i(n  - i)}{n^2} \big(n(n + 1) -
2i(n - i)\big)\,,
$$
and the coefficients $k_i$ in (5.1): $k_i = \dfrac{i(n  - i)}2$. It
follows that the ratio 
$$
\frac {\|\omega_i\|^2}{k_i^2} = \frac{4}{n^2} \left(\frac{n(n + 1)}{i(n  - i)} - 2\right)
$$
attains its maximum when $i = 1$ or $n-1$; therefore $k = 
\dsize{\frac{k_1}{\|\omega_1\|} = \frac{n}2\sqrt{\frac{n-1}{n^2 - n +
2}}}$. Similar computation can be done for root systems of other types.


\heading{\S 6. Geodesics and flats in locally symmetric spaces}\endheading 

\subhead{6.1}\endsubhead We are now going to use the result of the
previous section and  derive Theorems
1.4 and 1.10 from Theorems 1.7 and  1.9
respectively. Throughout the end of the section, $Y \cong K\backslash
\ggm$ is a   
noncompact irreducible 
locally symmetric space  of
noncompact type and finite volume. Here  
 $G$ is  the connected component of the identity in the isometry
group of the universal cover $\tilde Y$ of $Y$, 
$\Gamma$  an irreducible
lattice in $G$ and $K$ a maximal compact subgroup of $G$, i.e.~the
stabilizer of a point $\tilde y_0\in \tilde Y$. 
 
Denote by $\g$ (resp.~$\k$) the Lie algebra of $G$ (resp.~$K$). The
geodesic symmetry at $\tilde y_0$ induces a Cartan decomposition $\g =
\k \oplus \p$, and one can identify the tangent space to a point
$\tilde y_0\in Y$ with $\p$. Fix a Cartan
subalgebra  $\a$ of  $\p$. Let $\a_{\sssize +}$ be a positive Weyl chamber
relative to a fixed ordering of the root system of the pair $(\g,\a)$. 
Then it is known that the set $\a_1$ of unit vectors in $\a_{\sssize +}$ is a
fundamental set for the $G$-action on the unit tangent bundle
$S(\tilde Y)$ of $\tilde Y$; that is, every orbit of $G$ intersects
the set $\{(\tilde y_0,\vz)\mid \vz\in\a_1\}$ exactly
once. Furthermore, for $\vz\in\a_1$, the stabilizer of $(\tilde
y_0,\vz)$ in $G$ is the centralizer $K_\vz$ of $\vz$ in $K$, so the
$G$-orbit of $(\tilde
y_0,\vz)$ in $S(\tilde Y)$ (resp.~the
$G$-orbit $\Cal E_\vz \df G(y_0,\vz)$ of $(y_0,\vz)$ in $S(Y)$) can be
identified with $K_\vz\backslash G$ (resp.~with $K_\vz\backslash
\ggm$). The sets $\Cal E_\vz$, $\vz\in\a_1$, are smooth submanifolds
of $S(Y)$ of finite Riemannian volume, which form a singular
measurable foliation of  $S(Y)$. It will be convenient to introduce the
notation $\Cal E_{\vz,y}$ for the set of all  $\xi\in S_y(Y)$ for
which $(y,\xi)\in \Cal E_\vz$ (here $y$ is an arbitrary point of
$Y$). Note that if the $\br$-rank of $G$ is equal to $1$, the set $\a_1$
consists of a single element $\vz$, 
so one has $\Cal
E_\vz = S(Y)$ and $\Cal E_{\vz,y} = S_y(Y)$ for any $y\in Y$.

It has been shown by F.~Mautner \cite{Ma} that the geodesic flow
$\gamma_t$ on $S(Y)$
restricted to  $\Cal E_\vz$, $\vz\in\a_1$, can be described via the
action of the one-parameter subgroup $\{\exp(t\vz)\}$ of $G$ as
follows:
$$
\gamma_t(K_\vz g\Gamma) = K_\vz \exp(t\vz)g\Gamma\,.\tag 6.1
$$
This clearly provides a link between Theorems 1.4 and 1.7. In
particular, one can prove the following strengthening of Theorem 1.4:

\proclaim{6.2. Theorem} There exists ${k}  = {k} (Y)
> 0$ such that for any $\vz\in\a_1$ the following holds: if
$y_0\in Y$ and    $\{r_{t}\mid {t}\in
\bn\}$ is  a  sequence of real numbers, 
 then   for  any $y\in Y$ and almost
every (resp.~almost no) $\xi\in \Cal E_{\vz,y}$ there  are infinitely many
${t}\in \bn$
such that {\rm (1.2)} is satisfied, 
  provided the series $\sum_{{t}=1}^{\infty} e^{-{k} r_{t}}$
diverges (resp.~converges). \endproclaim

\demo{Proof} Let $p$ denote the natural projection from $X = \ggm$ onto
$\Cal E_\vz$, take $x_0\in p^{-1}(y_0)$ and denote by $\Delta$ the
function $\text{\rm
dist}_X(x_0,\cdot)$ on $X$. Using Proposition 5.1,
find $k$ such that  $\Delta$ is   $k$-DL.  If $\sum_{{t}=1}^{\infty}
e^{-{k} r_{t}} = \infty$,
then, by Theorem 1.7, for any $C > 0$ and almost all $x\in X$ there
are infinitely many ${t}\in \bn$ such that
$\Delta\big(\exp({{t}\vz})x\big) \ge  r_{t} + C$.
 But clearly $\Delta(x)$ and $\text{\rm
dist}_Y(y_0,y)$ differ by no more than  additive constant whenever
$p(x) = (y,\xi)$. Therefore it follows from (6.1) that the set
$$
\{(y,\xi)\in \Cal E_\vz\mid \text{ (1.2) holds for infinitely many
${t}\in\bn$}\}\tag 6.2
$$
has full measure in $\Cal E_\vz$. To
finish the proof of the divergence case, it remains to notice that 
for any $y,y' \in Y$ and $\xi\in \Cal E_{\vz,y}$  there exists
$\xi'\in \Cal E_{\vz,y'}$ such that $\text{\rm
dist}\big(\gamma_{t}(y,\xi), \gamma_{t}(y',\xi')\big)$ is uniformly
bounded from above for
all positive ${t}$. Therefore for any $y\in Y$ the intersection of the
set (6.2) with $\Cal E_{\vz,y}$ has full measure in the latter set.
The proof of the easier convergence case follows the same pattern
(and certainly it suffices to use Lemma 2.3 instead of the full
strength of Theorem 1.7).  \qed\enddemo

\demo{{\bf 6.3.}\ \ 
Proof of Theorem 1.4} The main statement is a direct consequence of
the above theorem and the decomposition of the volume measures on the
spheres $S_y(Y)$ in terms of the measures on the leaves  $\Cal
E_{\vz,y}$ for all $z\in \a_1$. As for 
 the logarithm law (1.3), its validity for the set  of pairs
$(y,\xi)$ of full measure in $S(Y)$ 
immediately follows  from Corollary 2.4 and
Lemma 2.8, and then, as in the above proof, one shows that the
intersection of this set with $S_y(Y)$ has full measure in $S_y(Y)$
for any $y\in 
Y$. 
 \qed\enddemo

\demo\nofrills{{\bf 6.4.}\ \ 
\it Proof of Theorem 1.10 }\ 
%
\ can be written along
the same lines, with minor modifications.   One considers the
$G$-action on the  bundle 
$S^d(\tilde Y)$ and finds a representative $(\vz_1,\dots,\vz_d)$, with
$\vz_i\in\a$, in any $G$-orbit (recall that $\p\supset\a$ has been identified
with the tangent space to $\tilde Y$ at $\tilde y_0$). Then $G$-orbits
in $S^d(Y)$ are identified with
quotients of $X = \ggm$ by centralizers $K_{(\vz_1,\dots,\vz_d)}$ in
$K$ of
appropriate ordered $d$-tuples $(\vz_1,\dots,\vz_d)$. Similarly
to (6.1), one describes 
$\gamma_\vt(K_{(\vz_1,\dots,\vz_d)} g\Gamma)$, where $\vt =
(t_1,\dots,t_d)\in \br^d$, via
the action of $\exp(\sum_i t_i 
\vz_i)$ on $K_{(\vz_1,\dots,\vz_d)}\backslash X$. An application of
Theorem 1.9 to the $\a$-action on $X$ then provides the desired
dichotomy, hence a logarithm law, for almost all 
$(y,\vec\xi)$ in any $G$-orbit. To derive a similar result for almost every
 $\vec\xi\in S_y^d(Y)$ and any $y\in Y$, one needs to
decompose $\a$ as a  union of Weyl chambers $\a_j$ and, accordingly,
break the flat  $\Cal F = \{\gamma_\vt(K_{(\vz_1,\dots,\vz_d)} g\Gamma)\mid \vt\in\d_{\sssize +}\} = \{K_{(\vz_1,\dots,\vz_d)} \exp(\sum_i t_i 
\vz_i)g\Gamma\mid \vt\in\d_{\sssize +}\}$ into 
pieces $\Cal F_j = \{K_{(\vz_1,\dots,\vz_d)} \exp(\sum_i t_i 
\vz_i)g\Gamma\mid \vt\in\d_{\sssize +},\, \sum_i t_i 
\vz_i \in \a_j\}$. After that it remains to notice that given each of
the pieces $\Cal F_j$ and a point 
$y\in Y$, one can use Iwasawa decomposition for $G$ to find a similar piece $\Cal F'_j$ starting from $y$
which lies at a bounded distance from $\Cal F_j$. \qed\enddemo

 \heading{\S 7. A very important DL function on the space of
lattices}\endheading

\subhead{7.1}\endsubhead We now consider another class of examples of
DL functions on homogeneous spaces. Throughout the section  we fix an
integer $k > 1$,  let $G = SL_k(\br)$, $\Gamma 
= SL_k(\bz)$ and $\mu$ the normalized Haar measure on the space $X_k \df
\ggm$ of unimodular lattices 
in $\br^k$, choose a norm on  $\br^k$ and define the function $\Delta$
on  $X_k$  by (1.9).
Our goal is to prove 

\proclaim{Proposition}  There exist positive  $C_k, C'_k$
such that 
$$
C_ke^{-kz} \ge \td(z) \ge C_ke^{-kz} - C'_k e^{-2kz} \quad\text{for all
} z \ge 0\,,\tag 7.1
$$
in particular, $\Delta$ is  $k$-DL.
 \endproclaim

The main tool here is the reduction theory for $SL_k(\br)/SL_k(\bz)$,
in particular, a 
generalization of Siegel's \cite{Si} summation formula. Recall that
a vector $\vv$ in a lattice $\Lambda\subset \br^k$ is called {\sl
primitive\/} (in $\Lambda$) if it is not a multiple of another element
of $\Lambda$; equivalently, if there exists a basis
$\{\vv_1,\dots,\vv_k\}$ of $\Lambda$ with $\vv_1 = \vv$. Denote by
$P(\Lambda)$ the set of all primitive vectors in $\Lambda$. Now, given
a function  $\ph$ on $\br^k$, define a function $\overset{\sssize \
\wedge}\to { \ph}$ on $X_k$ by $\overset{\sssize \wedge}\to {
\ph}(\Lambda)  \df 
\sum_{\vv\in P(\Lambda)}\ph(\vv)$. The following is one of the results of
the paper  \cite{Si}:

\proclaim{7.2. Theorem} For any  $\ph\in L^1(\br^k)$, one has
$
\int_{X_k}\overset{\sssize \wedge}\to { \ph}\,d\mu =
c_k\int_{\br^k}\ph\,d\vv\,,
$
where $c_k = \frac1{\zeta(k)}$.
 \endproclaim

The theorem below is a direct generalization of Siegel's result. 
For $1 \le {d} < k$, say that an ordered ${d}$-tuple $(\vv_1,\dots,\vv_{d})$
of vectors in a lattice 
$\Lambda\subset \br^k$ is  {\sl
primitive\/} if it is extendable to a basis of $\Lambda$, and denote by
$P^{{d}}(\Lambda)$ the set of all such  ${d}$-tuples. Now, given
a function  $\ph$ on $\br^{k{d}}$, define a function $\overset{\sssize
\ \wedge_{\sssize {d}}}\to { \ph\ \,}$ on 
$X_k$ by $\overset{\sssize \ \wedge_{\sssize {d}}}\to { \ph\ \,}(\Lambda) \df 
\sum_{(\vv_1,\dots,\vv_{d})\in
P^{{d}}(\Lambda)}\ph(\vv_1,\dots,\vv_{d})$. Then one has 

\proclaim{7.3. Theorem}  For $1 \le {d} < k$ and  $\ph\in L^1(\br^{k{d}})$, 
$$
\int_{X_k}\overset{\sssize \ \wedge_{\sssize {d}}}\to { \ph\ \,}\,d\mu =
c_{k,{d}}\int_{\br^{k{d}}}\ph\,d\vv_1\dots d\vv_{d}\,,\tag 7.2
$$ 
where $c_{k,{d}} = \frac1{\zeta(k)\cdot\dots\cdot\zeta(k-d+1)}\,$.
 \endproclaim

\demo{Sketch of proof} We essentially follow S.~Lang's
presentation ({Yale
University lecture course, Spring 1996}) of
Siegel's   
original proof. Fix a basis
$\{\ve_1,\dots,\ve_k\}$ of $\br^k$, denote by $G'$ (resp.~$\Gamma'$)
the stabilizer of 
the ordered ${d}$-tuple $(\ve_1,\dots,\ve_{d})$ in $G$ (resp.~in
~$\Gamma'$). Then $G/G'$, as  a $G$-\hs, can be naturally identified
with an open dense subset of $\br^{k{d}}$, namely, with the set of
 linearly independent $d$-tuples. Similarly $\Gamma/\Gamma'$ can be
identified with the 
$\Gamma'$-orbit of $(\ve_1,\dots,\ve_k)$, which is exactly the set
$P^{{d}}(\bz^k)$ of primitive  ${d}$-tuples in $\bz^k$.
These identifications allow one to transport the Lebesgue measure from
$\br^{k{d}}$ to a Haar measure $\mu_{G/G'}$ on $G/G'$, and to interpret
the summation over $P^{{d}}(\bz^k)$ as the integration over the counting
measure $\mu_{\Gamma/\Gamma'}$ on  $\Gamma/\Gamma'$.

The choice of the normalized Haar measure $\mu$ on $X_k$ (and hence
of the measures $\mu_{G}$ on $G$ and $\mu_{G/\Gamma'}$ on $G/\Gamma'$),
together with the aforementioned choice  of 
$\mu_{G/G'}$,  uniquely determines the Haar measures $\mu_{G'}$ and
$\mu_{G'/\Gamma'}$ on  $G'$ and $G'/\Gamma'$ (note that $\Gamma'$ is a
lattice in $G'$) such that for any $\ph\in
L^1(G/\Gamma')$ one has
$$
\int_{X_k}\int_{\Gamma/\Gamma'}\ph\,d\mu_{\Gamma/\Gamma'}\,d\mu =
\int_{G/\Gamma'}\ph\,d\mu_{G/\Gamma'} = 
\int_{G/G'}\int_{G'/\Gamma'}\ph\,d\mu_{G'/\Gamma'}\,d\mu_{G/G'}\,.\tag 7.3
$$

It remains to take any  $\ph\in L^1(\br^{k{d}})\cong L^1(G/G')$,
extend it to an 
integrable function on $G/\Gamma'$, and notice that the left hand side
of (7.2) coincides with that of (7.3), whereas the right hand side
of (7.3) can be rewritten as $\mu_{G'/\Gamma'}({G'/\Gamma'})\cdot
\int_{G/G'}\ph\,d\mu_{G/G'}$, which is exactly the  right hand side
of (7.2) with $c_{k,{d}} = \mu_{G'/\Gamma'}({G'/\Gamma'})$. The
computation of the exact value of $c_{k,{d}}$ is not needed for our
purposes and is left as an exercise for the reader. \qed\enddemo


\demo{{\bf 7.4.}\ \ 
Proof of Proposition 7.1} Take $z \ge 0$, denote by $B$ the ball in
$\br^k$ of radius $e^{-z}$ centered at the origin, and by $\ph$ the
characteristic function of $B$. Note that 
$$
\Delta(\Lambda) \ge z \thus \log\big(\tfrac1{\|\vv\|}\big)
\ge z \text{ for some }\vv\in\Lambda\nz \thus \Lambda\cap B \ne
\{0\}\,,
$$
and the latter condition clearly implies that $B$ contains at least
two primitive vectors ($\vv$ and $-\vv$) of $\Lambda$. Since
$\overset{\sssize \ 
\wedge}\to { \ph}(\Lambda) = \#\big(P(\Lambda)\cap B\big)$, one has  
$$
\int_{X_k}\overset{\sssize \wedge}\to { \ph}\,d\mu = \int_{\{\Lambda \mid
\Delta(\Lambda) \ge z\}} \overset{\sssize \wedge}\to { \ph}\,d\mu \ge 2\mu\big(\{\Lambda \mid
\Delta(\Lambda) \ge z\} \big)\,.\tag 7.5
$$
The left hand side, in view of  Theorem 7.2, 
is equal to $c_k\int_{\br^k}\ph\,d\vv = c_k\nu_k e^{-kz}$
(here $\nu_k$ is the volume of the unit ball in $\br^k$), hence the
upper estimate for $\td(z)$ in (7.1), with $C_k = \frac12  c_k\nu_k$. 

For the lower estimate, we will demonstrate that
lattices $\Lambda$ with $\overset{\sssize 
\wedge}\to { \ph}(\Lambda) > 2$ contribute very insignificantly to the
integral in the 
left hand side of 
(7.5). Indeed, a standard argument from reduction theory shows that
whenever there exist at least two linearly independent vectors in
$\Lambda\cap B$, for any $\vv_1\in P(\Lambda)$ one can find $\vv_2\in
\Lambda\cap B$ such that $(\vv_1,\vv_2)$, as well as $(\vv_1,-\vv_2)$,
belongs to
$P^2(\Lambda)$. Consequently, one has 
$$
\overset{\sssize \
\wedge}\to { \ph}(\Lambda) = \#\big(P(\Lambda)\cap B\big) \le
\frac12\#\big(P^2(\Lambda)\cap (B\times B)\big)
$$ 
whenever $\overset{\sssize \
\wedge}\to { \ph}(\Lambda) > 2$. Note that the right hand side is equal
to $\frac12\overset{\sssize \
\wedge_{\sssize 2}}\to { \psi}(\Lambda)$, where $\psi$ is the
characteristic function of $B\times B$ in $\br^{2k}$. Therefore
$$
\split
\int_{X_k}\overset{\sssize \wedge}\to { \ph}\,d\mu  &= \int_{\{\Lambda \mid
\overset{\sssize \
\wedge}\to { \ph}(\Lambda) = 2\}} \overset{\sssize \wedge}\to { \ph}\,d\mu   +  \int_{\{\Lambda \mid
\overset{\sssize \
\wedge}\to { \ph}(\Lambda) > 2\}} \overset{\sssize \wedge}\to { \ph}\,d\mu  \\
&\le 2\mu\big(\{\Lambda \mid
\overset{\sssize \
\wedge}\to { \ph}(\Lambda) = 2\}   \big) + \frac12\int_{\{\Lambda \mid
\overset{\sssize \
\wedge}\to { \ph}(\Lambda) > 2\}} \overset{\sssize \ \wedge_{\sssize
2}}\to { \psi}\,d\mu 
\le 2\mu\big(\{\Lambda \mid
\Delta(\Lambda) \ge z\} \big) + \frac12\int_{X_k} \overset{\sssize \ \wedge_{\sssize
2}}\to { \psi}\,d\mu \,.
\endsplit
$$
{}From Theorems 7.2 and 7.3 it then follows that  $2\td(z) \ge
c_k\nu_k e^{-kz} - \frac12c_{k,2}(\nu_k)^2 e^{-2kz}$, which finishes the proof of
the proposition. \qed
\enddemo

 

\heading{\S8. The Khinchin-Groshev Theorem}\endheading 

\subhead{8.1}\endsubhead We begin by introducing some terminology. Let
$\psi:\bn\mapsto (0,\infty)$ be a positive function.   Fix $m,n\in
\bn$ and say that a matrix $A\in\mr$  (viewed 
as a  system   of $m$ linear forms in $n$ variables) is {\sl
$\psi$-approximable\footnote{The authors are grateful to M.~Dodson for
 a permission to modify his terminology introduced in
\cite{Do}.}\/} if  there are infinitely many $\vq\in \bz^n$ such that
(1.1) holds. 
Then one can restate Theorem
1.1 as follows:

\proclaim{Theorem}  Let $m$,$n$ be positive integers and
$\psi:[1,\infty)\mapsto (0,\infty)$ a 
non-increasing continuous function. Then almost every (resp.~almost no)
$\Lambda\in X_\mn$  is 
$(\psi,n)$-approximable,  provided the integral $
{\int_{1}^\infty {\psi(x)}\,dx} 
$
 diverges (resp.~converges).
\endproclaim  

To prepare for the reduction of this theorem to Theorem 1.7,  let us
present an equivalent formulation.  For a vector $\vv\in \br^\mn$,
denote by  $\vv^{(m)}$ 
(resp.~$\vv_{(n)}$) the 
vector consisting of first $m$ (resp.~last $n$)  components of $\vv$. 
Now, to  a matrix
$A\in\mr$ we associate a lattice $\la$ in $\br^\mn$ defined by $\la
\df \left(\matrix 
I_m & A  \\
0 & I_n
\endmatrix \right)\bz^\mn$; in other words, $\la = \left\{\left.\left(\matrix
A\vq + \vp  \\
\vq
\endmatrix \right)\right|\vp\in\bz^m,\vq\in\bz^n\right\}$. 
Clearly   $A\in\mr$ is $\psi$-approximable iff  there exist
$\vv\in \la$ with  arbitrarily large $\|\vv_{(n)}\|$ such that 
$$
 \|\vv^{(m)}\|^m   \le \psi(\|\vv_{(n)}\|^n)\,.\tag 8.1
$$

Let us say that a lattice $\Lambda \in X_\mn$ is  {\sl $(\psi,n)$-approximable}
iff  there exist $\vv\in \la$  with   arbitrarily large
$\|\vv_{(n)}\|$  such that (8.1) 
holds. Now the above theorem   can be restated as follows: 

\roster 
\item"{$\bullet$}"\it  Let $m$,$n$ be positive integers and
$\psi:[1,\infty)\mapsto (0,\infty)$ a 
non-increasing continuous function. Then 
almost every 
(resp.~almost no)  lattice of the form $\la$, $A\in\mr$,  is
$(\psi,n)$-approximable, provided the integral 
$
{\int_{1}^\infty {\psi(x)}\,dx} 
$ diverges (resp.~converges).
\endroster 

\rm
We will see
later that the 
same phenomenon takes place  for generic lattices in $\br^\mn$. More
precisely, we will prove 

\proclaim{8.2. Theorem}  Let  $\psi$, $m$ and $n$ be as in Theorem 8.1. 
Then almost every (resp.~almost no) $\Lambda\in X_\mn$  is
$(\psi,n)$-approximable,  provided the integral $
{\int_{1}^\infty {\psi(x)}\,dx} 
$
 diverges (resp.~converges).\endproclaim  

In fact it is not a priori clear how to derive Theorem 8.2 from Theorem
1.1 and vice versa. We will do it 
by restating
these theorems in the language of flows on the space of lattices.
For that we first 
need a change of variables 
technique formalized in the following 
 
\proclaim{8.3. Lemma}  Fix $m,n\in\bn$  and  $x_0 > 0$, and let
$\psi:[x_0,\infty)\mapsto (0,\infty)$ be a non-increasing continuous
function. Then there exists a unique 
continuous  function  $r:[t_0,\infty)\mapsto \br$, where
$t_0 =  \frac m\mn \log x_0 -\frac n\mn \log \psi(x_0)$, such that  
$$
\text{the function} \quad {\lambda(t)\df t-nr(t)}  \quad \text{is strictly increasing and tends to $\infty$ as } t\to+\infty\,,
\tag 8.2a
$$
$$
\text{the function} \quad {L(t) \df t+mr(t)}  \quad \text{is nondecreasing}\,,
\tag 8.2b
$$
and
$$
\psi(e^{t-nr(t)}) =  e^{-t - m r(t)}\quad\forall\,t\ge t_0\,.\tag 8.3
$$
Conversely, given $t_0 \in\br$ and a continuous  function
$r:[t_0,\infty)\mapsto\br$ such that {\rm (8.2ab)} hold, there exists
a unique continuous non-increasing function
$\psi:[x_0,\infty)\mapsto(0,\infty)$, with $x_0 = e^{t_0 - nr(t_0)}$,
satisfying {\rm (8.3)}. 
Furthermore, for a nonnegative integer $q$, 
$$
I_1 \df \int_{x_0}^\infty {(\log x)^{q}\psi(x)}\,dx < \infty\quad
\text{iff} \quad
I_2 \df \int_{t_0}^\infty t^{q} e^{-(m+n)r(t)}\,dt < \infty\,.
$$ \endproclaim

\demo{Proof} The claimed correspondence becomes transparent if one
uses the variables \linebreak $L = -\log \psi$, $\lambda = \log
x$, and the function $P(\lambda) \df -\log \psi(e^\lambda)$
. Given $t\ge t_0$, one can define $\big(\lambda(t),L(t)\big)$ to be 
the unique intersection point of the graph of the nondecreasing
function $L = P(\lambda)$ and the decreasing straight line
$L = \frac\mn n t - \frac mn\lambda$, and then put
$$
r(t) = \frac{L(t)
- \lambda(t)}\mn\,.\tag 8.4
$$
 The properties (8.2ab) and (8.3) are then
straightforward.  Conversely, given the function $r(\cdot)$ with
(8.2ab) and $\lambda \ge \lambda_0\df t_0 - nr(t_0)$, one defines $P(\lambda)$ to
be equal to $L\big(t(\lambda)\big)$, where $L(\cdot)$ is as in
(8.2b) and  $t(\cdot)$ is the function inverse to  $\lambda(\cdot)$ of
(8.2a).  

Further, the integral $I_1$ is equal to $\int_{\lambda_0}^\infty
\lambda^{q}e^{\lambda-P(\lambda)}\,d\lambda$, while $I_2$, in view of
(8.2ab) and (8.4),  can be written as
$
\int_{\lambda_0}^\infty  \big(\frac m\mn\lambda + \frac n\mn
P(\lambda)\big)^{q}e^{\lambda-P(\lambda)}\,\big(\frac m\mn d\lambda +
\frac n\mn dP(\lambda)\big) \ge I_1$.
 It remains to assume $I_1 < \infty$ and prove that the following
integrals are finite: 
$$
I_3 = \int_{\lambda_0}^\infty  \lambda^{q}e^{\lambda-P(\lambda)}\,
dP(\lambda), \quad 
 I_4 = \int_{\lambda_0}^\infty
P(\lambda)^{q}e^{\lambda-P(\lambda)}\,d\lambda, \quad 
 I_5 = \int_{\lambda_0}^\infty
P(\lambda)^{q}e^{\lambda-P(\lambda)}\, dP(\lambda)\,. 
$$
Integration by parts reduces $I_3$ to  the form
$$
I_3 = - \int_{\lambda_0}^\infty  \lambda^{q}e^{\lambda}\,d\big(e^{-P(\lambda)}\big) = -\left.
\lambda^{q}e^{\lambda-P(\lambda)}\right|_{\lambda_0}^\infty +
\int_{\lambda_0}^\infty  e^{\lambda}(\lambda^{{q}} +
q\lambda^{{q}-1})e^{\lambda-P(\lambda)}\,d\lambda\,, 
$$
where both terms are finite due to the finiteness of $I_1$. To
estimate $I_4$, one writes 
$$
 I_4 = \int_{\lambda\ge \lambda_0,\,P(\lambda) < 2\lambda}
P(\lambda)^{q}e^{\lambda-P(\lambda)}\,d\lambda +  \int_{\lambda\ge
\lambda_0,\,P(\lambda) \ge 
2\lambda}   P(\lambda)^{q}e^{\lambda-P(\lambda)}\,d\lambda\,; 
$$
the first term is clearly bounded from above by $2^qI_1$, while the
integrand in the second term is for large enough values of $\lambda$
not greater than $2^q\lambda^{q}e^{-\lambda}$. This implies that $I_4$
is also finite. Finally, 
$$
\split
 I_5 &= \int_{\lambda\ge \lambda_0,\,P(\lambda) < 2\lambda}
P(\lambda)^{q}e^{\lambda-P(\lambda)}\,dP(\lambda) +  \int_{\lambda\ge
\lambda_0,\,P(\lambda) \ge 
2\lambda}   P(\lambda)^{q}e^{\lambda-P(\lambda)}\,dP(\lambda) \\ &\le 2^q
I_3 + \int_{\lambda_0}^\infty
P(\lambda)^{q}e^{-P(\lambda)/2}\,dP(\lambda) <\infty\,, 
\endsplit
$$
 which finishes the proof of the lemma. \qed
\enddemo

In what follows, we will denote by $\Cal D_{m,n}(\psi)$ (after S.G.~Dani) the
function $r$ 
corresponding to $\psi$  by the above lemma.  Note also that
$r$ does not have to be monotonic, but is always
quasi-increasing (as defined in \S 2.4) in view of (8.2b). 

\example{8.4. Example} The easiest special case is given by $\psi(x) =
\vre/x$ for a positive constant $\vre$. Then the equation (8.3) gives
$r(t) = \frac1{\mn}\log(\frac1\vre)$, so the correspondence $\Cal D_{m,n}$
sends such a function $\psi$ to a constant. Recall that \amr\ is said
to be {\sl badly 
approximable\/} if it is not $\frac\vre x$-approximable  for some $\vre >
0$. In \cite{D}, Dani proved that $A$ is \ba\ iff the trajectory
$\{f_t \la\mid t \ge 0\}$, with $f_t$ as in (1.10), is bounded in
$X_\mn$. Note that in view of Mahler's Compactness Criterion (see
\cite{R1, Corollary 10.9}), the latter condition is equivalent to the existence of an
upper bound for $\{\Delta(f_t\la) \mid t \ge 0\}$, with $\Delta$ as in
(1.9). 
\endexample

\subhead{8.5}\endsubhead We are now going to prove a generalization of
the aforementioned  result of Dani.

\proclaim{Theorem}  Let  $\psi$, $m$ and  $n$  be as in Theorem 8.1,   $\Delta$ as in
{\rm (1.9)}, $\{f_t\}$ as in {\rm (1.10)} Then $\Lambda\in X_\mn$ is $(\psi,n)$-approximable iff   there exist arbitrarily large positive  $t$
such that 
$$
\Delta(f_{t}\Lambda) \ge r(t)\,,\tag 8.5
$$ 
where
$r = \Cal D_{m,n}(\psi)$. In particular, \amr\  is
$\psi$-approximable iff    there exist arbitrarily large
positive  $t$ 
such that  
$$
\Delta(f_{t}\la) \ge r(t)\,.\tag 8.5$A$
$$
 \endproclaim

\demo{Proof} Assume that (8.1) holds for some $\vv\in\Lambda$, and
note that, by definition of 
$f_t$ and $\Delta$, to prove (8.5) it
suffices to find $t$ such  that 
$$
e^{t/m} \|\vv^{(m)}\| \le e^{-r(t)}\tag 8.6a
$$
and
$$
e^{-t/n} \|\vv_{(n)}\| \le e^{-r(t)}\tag 8.6b
$$
Now define 
$t$ by 
$$
\|\vv_{(n)}\|^n =
e^{t-nr(t)}\,.\tag 8.7
$$
In view of (8.2a), one can do this   whenever $
\|\vv_{(n)}\|$ is large enough. 
 Then (8.6b) follows immediately, and one can write
$$
 \|\vv^{(m)}\|^m   \un{(8.1)}{\le}
\psi(\|\vv_{(n)}\|^n)  \un{(8.7)}{=} \psi(e^{t-nr(t)}) \un{(8.3)}{=}
e^{-t - mr(t)}\,,
$$
which readily implies (8.6a). Lastly, again in  view of (8.2a), $t$
will be  arbitrarily large if one chooses $
\|\vv_{(n)}\|$  arbitrarily large as well.

For the converse, let us first take care of the case when
$$
\vv^{(m)} = 0\text{ for some }\vv\in\Lambda\nz\,.\tag 8.8
$$
 Then one can take
integral multiples of this $\vv$ to produce infinitely many vectors
satisfying (8.1); thus  lattices with (8.8) are $(\psi,n)$-approximable
for any function $\psi$. Otherwise, assume that  (8.5) holds for a
sufficiently large 
$t$. This immediately gives  a vector $\vv 
\in\Lambda$ satisfying (8.6a) and (8.6b), and one can write 
$$
 \|\vv^{(m)}\|^m  \un{(8.6a)}{\le}
e^{-t - mr(t)}\un{(8.3)}{=} \psi(e^{t-nr(t)}) \un{(8.6b) and the
monotonicity of $\psi$}{\le} \psi(\|\vv_{(n)}\|^n)\,.
$$ 
Finally, if $t$ is taken arbitrarily large, $\|\vv^{(m)}\|$ becomes
arbitrarily small in view of (8.6a), and yet can not equal  zero, so
$\|\vv_{(n)}\|$ must be  
arbitrarily large  by the discreteness of $\Lambda$.  \qed\enddemo

\demo{{\bf 8.6.}\ \ 
Proof of Theorem 8.2} In view of the above theorem and Lemma 8.3, it
suffices to  prove the
following

\proclaim{Theorem} Given $m,n\in \bn$,  $\Delta$ as in
{\rm (1.9)}, $\{f_t\}$ as in {\rm (1.10)} and a continuous
quasi-increasing function $r:[t_0,\infty)\mapsto\br$, for 
almost every (resp.~almost no) $\Lambda\in X_\mn$  there exist
arbitrarily large positive  $t$ 
such that {\rm (8.5)} holds,  provided the integral 
$
\int_{t_0}^\infty  e^{-(m+n)r(t)}\,dt
$
 diverges (resp.~converges). 
 \endproclaim

\demo{Proof} From Corollary 2.4 and Lemma 2.8 it is clear that the
above statement 
is a straightforward  consequence of the family  $\cb(\Delta)$
being Borel-Cantelli for $f_1$. The latter, in its turn, immediately follows
 from Theorem 1.7 and Proposition 7.1. \qed\enddemo
\enddemo

\demo{{\bf 8.7.}\ \ 
Proof of Theorem 1.1}  Similarly, Theorem 1.1 follows from

\proclaim{Theorem} 
Given $m,n\in \bn$,  $\Delta$ as in
{\rm (1.9)}, $\{f_t\}$ as in {\rm (1.10)} and a continuous
quasi-increasing function $r:[t_0,\infty)\mapsto\br$, for 
almost every (resp.~almost no) \amr\  there exist
arbitrarily large positive  $t$ 
such that {\rm (8.5$A$)} holds,  provided the integral $
\int_{t_0}^\infty  e^{-(m+n)r(t)}\,dt
$
 diverges (resp.~converges). \endproclaim

\demo{Proof} 
It is easy to see (cf.~\cite{D, 2.11}) that any lattice $\Lambda\in
X_\mn$ can be written in the form  
$$
\Lambda = \left(\matrix
B_1 & 0  \\
B_2 & B_3
\endmatrix \right)\la\,,
$$
for some \amr, $B_1\in M_{m,m}(\br)$, 
$B_2\in M_{n,m}(\br)$ and $B_3\in M_{n,n}(\br)$ with \linebreak 
$\text{det}(B_1)\text{det}(B_3) = 1$. Therefore one can write
$$
f_t\Lambda = f_t\left(\matrix
B_1 & 0  \\
B_2 & B_3
\endmatrix \right)f_{-t}f_t\la = \left(\matrix
B_1 & 0  \\
e^{-(t/m + t/n)}B_2 & B_3
\endmatrix \right)f_t\la\,.
$$
{}From this and the uniform continuity of $\Delta$ it follows that for
some positive $C$ (dependent on $\Lambda$) one has $\sup_{t>
0}|\Delta(f_t\Lambda) - \Delta(f_t\la)| < C$. If 
$\int_{t_0}^\infty  e^{-(m+n)r(t)}\,dt
$
diverges (resp.~converges), by Theorem 8.6 the set of lattices  
$\Lambda$ such that for any (resp.~for some) $C > 0$   
there exist
arbitrarily large positive  $t$ with $
\Delta(f_{t}\Lambda) \ge r(t) + C$ (resp.~with $
\Delta(f_{t}\Lambda) \ge r(t) - C$), has full (resp.~zero) measure in
$X_\mn$. Therefore, by Fubini, the set of \amr\ such that (8.5$A$) holds
for arbitrarily large $t$ has full (resp.~zero) measure in  $\mr$. \qed\enddemo
\enddemo

\example{8.8. Remark} It is also possible to argue in the opposite
direction and deduce Theorem 8.6
from Theorem 8.7. (Cf.~\cite{D}, where the abundance of bounded orbits for
certain flows on $X_\mn$ was deduced from W.~Schmidt's result on \ba\
systems of linear forms, vs.~\cite{KM}, where ergodic theory was used
to construct bounded orbits, thus providing another proof of the
aforementioned result of Schmidt.) In
other words, one can derive logarithm laws for specific flows on
$X_\mn$ simply by applying Theorem 8.5 to translate the
Khinchin-Groshev Theorem into the 
dynamical language. As a historical note, the authors want to point
out that this is exactly what they understood first and what prompted
them to start working on this circle of problems. 
\endexample

\heading{9. Multiplicative approximation of lattices}\endheading

\subhead{9.1}\endsubhead As a motivation, let us consider the case $m
= n = 1$ of Theorem 8.2. The 
inequality (8.1) then transforms into  
$$
|v_1| \le \psi(|v_2|)\,,\quad \text{or} \quad |v_1||v_2| \le
|v_2|\psi(|v_2|)\,,\tag 9.1 
$$
where $\vv = (v_1,v_2)$ is a vector from a lattice $\Lambda\in X_2$.
 Since $\psi$ is bounded, the fact that (9.1) holds for vectors $\vv$
with arbitrarily large $|v_2|$ implies that one has $\|\vv\| = |v_2|$
for infinitely many $\vv\in\Lambda$ satisfying (9.1); therefore (9.1)
can be replaced by (1.11).  Conversely, if 
(1.11)  holds for  infinitely many $\vv\in\Lambda$, then either $\Lambda$
or  its mirror reflection around the axis $v_1 = v_2$ is
$(\psi,1)$-approximable.
This way one gets an equivalent form  of the
$m = n = 1$  case  of Theorem 8.2 as follows:  
  
\roster 
\item"{$\bullet$}" {\it  With $\psi$ as in Theorem 8.1, for almost every
(resp.~almost no)  $\Lambda\in X_2$  
there exist infinitely many $\vv\in
\Lambda$ with {\rm (1.11)},
   provided the integral $
{\int_{1}^\infty {\psi(x)}\,dx} 
$ diverges (resp.~converges).}
\endroster

This suggests a natural generalization and (sigh!) another definition:
for  an integer $k \ge 2$, say that $\Lambda\in X_k$ is {\sl
$\psi$-multiplicatively approximable\/} (to be abbreviated as $\psi$-MA)
if   there exist infinitely many $\vv\in \Lambda$ satisfying (1.11).
 Thus the above theorem can be restated as follows:

\roster 
\item"{$\bullet$}" {\it   For $\psi$ as in Theorem 8.1,  almost every
(resp.~almost no)  $\Lambda\in X_2$  is $\psi$-MA,   provided the integral 
 $
{\int_{1}^\infty {\psi(x)}\,dx} 
$  diverges (resp.~converges).} 
\endroster

 A question, raised by M.~Skriganov in \cite{Sk, p.~23}, amounts to
considering a family of functions 
$\psi_q(x) = 1/{x(\log x)^{q}}$ and looking
for a critical exponent $q_0 = q_0(k)$ such that almost all
(resp.~almost no)  $\Lambda\in
X_k$  are $\psi_q$-MA if $q \le q_0$ (resp.~if $q > q_0$). It is shown
in \cite{Sk} that $q_0(k)$ 
must be positive and not greater than $k - 1$. In this section we 
prove Theorem 1.11, which, using the above terminology,  reads
as follows:

\roster 
\item"{$\bullet$}" {\it   Let $\psi:[1,\infty)\mapsto (0,\infty)$ be a 
non-increasing continuous function and $k$ an integer greater than $1$. Then
almost every 
(resp.~almost no)  $\Lambda\in X_k$  is $\psi$-MA,   provided the integral 
$
{\int_{1}^\infty {(\log x)^{k -2}\psi(x)}\,dx} 
$ diverges (resp.~converges).}  
\endroster 

 In particular,  this proves the existence of  $q_0(k)$ and gives its
exact value, namely, $q_0(k)  = k - 1$.

\subhead{9.2}\endsubhead In order to reduce Theorem 1.11 to Theorem
1.9, we need an analogue of 
the correspondence of Theorem 8.5. Again, the special case given by
$\psi(x) = \vre/x$ and $r \equiv  \text{const}$  is worth mentioning.  Recall
that $\Lambda$ is called 
{\sl admissible\/} (cf.~\cite{Sk, p.~6}) if it is not $\frac\vre x$-MA  for
some $\vre > 
0$.   It easily
follows from Mahler's Compactness Criterion  (and is 
mentioned in \cite{Sk, p.~14}) that a lattice is 
admissible iff its orbit under the diagonal subgroup of $SL_k(\br)$ is
bounded in $X_k$. To generalize this observation, identify the Lie
algebra $\d$ of traceless diagonal $k\times k$ matrices  with $\{\vt =
(t_1,\dots,t_k)\in 
\br^k\mid\sum_{i = 1}^k t_i = 0\}$, denote by $f_\vt$ the element
of $SL_k(\br)$ given by 
$$
f_\vt = \exp(\vt) = 
\text{diag}(e^{t_1},\dots,e^{t_k})\,,\tag 9.2
$$
 and let $\|\vt\|_{\sssize -} \df  \max\{|t_i|\bigm|t_i
\le 0\}$.

\proclaim{Theorem}  Let  $\psi$ be as in Theorem 8.1, $k$ an
integer greater than $1$, 
$\Delta$ as in  
{\rm (1.9)} and $\{f_\vt\}$ as in {\rm (9.2)}. Then 
$\Lambda\in X_k$  is
$\psi$-MA   iff    there exist
$\vt\in\d$ arbitrarily 
far from $0$
such that 
$$
\Delta(f_{\vt}\Lambda) \ge r(\|\vt\|_{\sssize -})\,,\tag 9.3
$$ 
where
$r = \Cal D_{k-1,1}(\psi)$.  \endproclaim

\demo{Proof}  Assume that (1.11) holds for some $\vv\in\Lambda$; our
goal is to find $\vt$ such that 
$$
e^{t_i} |v_i| \le e^{-r(\|\vt\|_{\sssize -})} \text{ for all } 1 \le i \le k\,.\tag 9.4
$$ 
We will do it in two steps. First define $t\in\br$ by $\|\vv\| =
e^{t-r(t)}$ (as before, one uses  (8.2a) to justify this step  if $
\|\vv\|$ is large enough). Note that in view of (8.3) one then has
$$
\psi(\|\vv\|)  = \psi(e^{t-r(t)}) =  e^{-t-(k-1)r(t)}\,.
$$
To define $\vt$, assume without loss of generality that $|v_i| \ge
|v_{i+1}|$ for all $i < k$,  and put
$
e^{t_1} = \frac {e^{-r(t)}}{|v_1|} = \frac {e^{-r(t)}}{\|\vv\|} =
e^{-t}\,,
$
and then, inductively, 
$
e^{t_i} = \min\big(\frac {e^{-r(t)}}{|v_i|}, e^{-(t_1 + \dots +
t_{i-1})}\big)\,.
$
Then one can check that: 

\roster  
\item"{$\bullet$}" $
e^{t_i}$ is not greater than $\frac {e^{-r(t)}}{|v_i|}$ for all $i$,

\item"{$\bullet$}" $\sum_{i = 1}^k t_i = 0$, and

\item"{$\bullet$}"  
$t = -t_1 = - \min_{ 1 \le i \le k}t_i = \|\vt\|_{\sssize -}$\ . 
\endroster

 Therefore (9.4)
is satisfied, and it remains to observe that,  again in  view of
(8.2a), $\|\vt\|_{\sssize -}$
will be  arbitrarily large if one chooses $
\|\vv\|$  arbitrarily large as well.

For the converse, we have to first take care of the case when
$$
\vv_i = 0\text{ for some }\vv\in\Lambda\nz \text{ and }1 \le i \le k\tag 9.5
$$
(in \cite{Sk} such lattices are called {\sl not weakly
admissible\/}). Clearly one can take 
integral multiples of this $\vv$ to produce infinitely many vectors
satisfying (1.11); thus  lattices with (9.5) are $\psi$-MA
for any function $\psi$. Otherwise, assume that  (9.3) holds for some
$\vt\in\d$  with
sufficiently large 
$\|\vt\|_{\sssize -}$. This immediately gives  a vector $\vv 
\in\Lambda$ satisfying (9.4). Let us again order the
components of $\vv$ so that $|v_1| \ge \dots \ge 
|v_{k}|$. Note that  without loss of generality one can assume that
$\|\vt\|_{\sssize -} = -t_1$ (otherwise, if $\|\vt\|_{\sssize -} =
-t_j > -t_1$,  one can interchange $t_1$ and $t_j$ without any damage
to (9.4)).  Now one can  multiply the
inequalities (9.4) for $i = 2, \dots,n$ by each other  to get 
$
\prod_{2 \le i \le k}e^{t_i} |v_i| \le e^{-(k-1)r(\|\vt\|_{\sssize
-})}$, which makes $
  {\Pi(\vv)}/{ \|\vv\|}$ to be not greater than
$$
  e^{t_1-(k-1)r(\|\vt\|_{\sssize -})} = e^{-\|\vt\|_{\sssize
-}-(k-1)r(\|\vt\|_{\sssize -})} 
\un{(8.3)}{=}  \psi(e^{\|\vt\|_{\sssize -} -r(\|\vt\|_{\sssize
-})})\un{(9.4) and the 
monotonicity of $\psi$}{\le}  \psi(\|\vv\|)
$$   
as desired.
Finally, recall that $\vt$ can be taken arbitrarily far from $0$. Let
$i$ be such that $t_i = \max_{1\le j \le k} t_j$.  Then (9.4) makes
$|v_i|$ 
arbitrarily small and yet positive, so
$\|\vv\|$ must be  
arbitrarily large  by the discreteness of $\Lambda$.  \qed\enddemo

\demo{{\bf 9.3.}\ \ 
Proof of Theorem 1.11} In view of the correspondence described in the
above theorem, we have to prove the
following

\roster 
\item"{$\bullet$}"   {\it Given an integer $k > 1$, $\Delta$ as in
{\rm (1.9)},  $\d$ as in \S 9.2, $\{f_\vt\}$ as in {\rm (9.2)} and a
continuous 
quasi-increasing function $r:[t_0,\infty)\mapsto\br$, for 
almost every (resp.~almost no) $\Lambda\in X_k$  there exist $\vt\in\d$ arbitrarily
far from $0$
such that {\rm (9.3)} holds,  provided the integral 
$
\int_{t_0}^\infty t^{k-2} e^{-kr(t)}\,dt 
$
 diverges (resp.~converges).}
\endroster 

Note that  the
function $\vt\mapsto\|\vt\|_{\sssize -}$ becomes a norm when restricted to any Weyl
chamber of $\d$.  Therefore  one can
decompose $\d$ as a union of such chambers $\d_j$ and then apply
Theorem 1.9, powered by Proposition 7.1, to  conclude that  the
family $\cb(\Delta)$ is Borel-Cantelli for 
$\{f_\vt\}$, where $\vt$ runs through the intersection of $\d_j$ with
an arbitrary 
lattice in $\d$. The statement of the theorem then immediately follows
 from  Corollary 2.4 and Lemma 2.8. \qed\enddemo 


\heading{\S 10. Concluding remarks and open questions}\endheading

\subhead{10.1}\endsubhead  It seems natural to conjecture that
the conclusion of Theorem 1.12 (isolation properties of the
restriction   of $\rho_0$ to any simple  
factor of $G$), and hence of Corollary 3.5 (exponential decay of
correlation coefficients of smooth functions),   
are satisfied for uniform lattices $\Gamma\subset G$ as well. This is
clearly the case when all factors of $G$ have property (T); otherwise
the problem stands  open. 

\subhead{10.2}\endsubhead  In view of the result of W.~Philipp
mentioned in \S 1.5 (or a similar result
for expanding rational maps of Julia sets announced recently by
R.~Hill and S.~Velani), it seems
natural to ask whether the family  of all balls in $\ggm$ will be
Borel-Cantelli for an element $\exp(\vz)$  of $G$ as in Theorem 
1.7.  For fixed $x_0\in\ggm$, this would measure the rate with which
a typical orbit approaches $x_0$, in particular, in the form of a
\lol\ for the function 
$
\Delta (x) = \log\big(\frac1{\text{dist}(x_0,x)}\big)
$.
This function satisfies ($k$-DL) with $k = \text{dim}(\ggm)$, but
is not uniformly continuous, therefore super-level sets of $\Delta$
cannot be adequately approximated by smooth functions. 

On the other
hand, D.~Dolgopyat \cite{Dol} recently proved a number of limit
theorems for partially hyperbolic dynamical systems. In particular he
showed that if $f$ is a partially hyperbolic diffeomorphism of a
compact Riemannian 
manifold $X$, then   the family  of
all balls in $X$ 
is
Borel-Cantelli for $f$, provided a certain additional assumption 
(involving rate of convergence of averages along pieces of unstable
leaves) is satisfied. 
Using  \cite{KM, Propositions 2.4.8 or A.6} this assumption can be
checked when $G$, $\Gamma$  and $f = \exp(\vz)$ are  as in Theorem
1.7,  $X = \ggm$ is compact and all simple factors of $G$ have
property (T). See also 
\cite{CK, CR} for other results in 
this direction.

\subhead{10.3}\endsubhead  We now roughly sketch modifications one
has to make in order to consider flows on reducible \hs s. If $G$
is a  connected semisimple center-free 
Lie 
group  without compact factors and  $\Gamma$ is a
lattice in $G$, one can find  connected normal subgroups
$G_1,\dots,G_l$ of $G$ such 
that 
$G=\prod_{i=1}^{l}G_i$ (direct product), 
$\Gamma_i \df G_i\cap\Gamma$ is an irreducible lattice in $G_i$ for each
$i$,  and $\prod_{i=1}^l\Gamma_i$ has finite index in $\Gamma$
(cf.~\cite{R1, Theorem 5.22}). As
a consequence of the above, 
$\ggm$ is finitely covered by  the direct product of  the spaces
$G_i/\Gamma_i$.  Denote by  $p_i$
the projection from $G$ onto  $G_i$. Then one can apply Corollary	
3.5 to the factors $G_i/\Gamma_i$ (more precisely, to the		
noncompact ones) and deduce that Theorem 4.3 (hence Theorem 1.8
as well) holds in
this generality provided the condition (ED) is replaced by 
$$
p_i(F)\text{ is  
{ED}  for all }i = 1,\dots,l\,.\tag 10.1
$$
Similarly one takes care of the case when $G$ has a nontrivial
center: then $G$ can be written as an almost direct product of the groups
$G_i$, and the maps $p_i$ are defined to be the projections $G\mapsto
G/\prod_{j\ne i}G_j$. 

Specializing to the case $F = \{\exp({{t}\vz})\mid t\in\bn\}$, with
$\vz\in \a$ as in Theorem 1.7, it is easy to see that (10.1) is
satisfied whenever $p_i(\vz)$ is nontrivial for all $i$ (here with
some abuse of notation we let $p_i$ be the projections of the
corresponding Lie algebras). The latter condition 
  holds for a
generic element $\vz\in\a$. Furthermore, one can prove that the
$k$-DL property of the distance function can be lifted to the
direct product of metric spaces.  (More precisely, if $\Delta_i$ is a
$k_i$-DL function on $X_i$, $1\le i \le l$, then $\sqrt{\Delta_1^2
+ \dots +\Delta_l^2}$ is $(\min_{1\le i \le l}k_i)$-DL function on
$\prod_{i=1}^{l}X_i$.) Therefore one can argue as in \S 6 and
prove  Theorem 1.4  without assuming that the space $Y$ is irreducible.

\subhead{10.4}\endsubhead Suppose that $G$,  $\Gamma$ and  $F =
\{f_{t}\}$ are as in Theorem 1.8, and let $\Delta$ be a uniformly
continuous function on $\ggm$ such that 
$$
\forall\, c < 1 \ \exists\,\delta > 0 
\text{ such that }\td(z + \delta) \ge c\cdot\td(z) \text{ for large
enough }z\,.\tag 10.2
$$
For such functions one can prove a refinement of Theorem 4.3:  if
$\{r_{t}\}$ is a sequence of real numbers satisfying (4.5), then
for  almost 
all $x\in \ggm$ one has
$$
 \lim_{N\to\infty}\frac {\#\{ 1 \le {t}\le N\mid \Delta(f_{t}x)
\ge r_{t}\}}  {\sum_{s=1}^{N} \td\big(r_{t}\big)} = 1\,.
$$
It is easy to see that (7.1) implies (10.2), therefore such a
refinement holds for the function $\Delta$ on
$SL_{k}(\br)/SL_{k}(\bz)$ given by (1.9). It seems very likely
that distance functions on locally symmetric spaces satisfy (10.2)
as well; in other words, one can write  exact asymptotics for the
measure of the complement of a ball of radius $z$, not only 
bound it from both sides by const$\cdot e^{-kz}$. However, the
proof is beyond our reach, since in order to use the main tools of
our proof (reduction
theory and the quasi-isometry with a Siegel set) one has to
sacrifice a multiplicative constant.

\heading 
{Appendix}
\endheading

\subhead{A.0}\endsubhead Let $\rho$ be a  unitary representation of a  locally compact second countable
group $G$ in a separable Hilbert space $V$. Say that a
sequence $\{v_{t}\mid t\in\bn\}\subset V$ 
is {\sl asymptotically
$\rho$-invariant\/} if $v_{t} \ne 0$ for all sufficiently large ${t}$, and $\|\rho(g)v_{t} - v_{t}\|/\| v_{t}\|\to 0$ as ${t}\to \infty$
uniformly 
on compact subsets of $G$. Then $\rho$ is isolated from $I_G$
 in the Fell topology iff there are no  asymptotically
$\rho$-invariant sequences $\{v_{t}\}\subset V$. 

Let now  $(X,\mu)$ be a probability space, and $(g,x)\mapsto gx$ a
$\mu$-preserving action of $G$ on $X$. Denote by $L^2_0(X,\mu)$ the subspace  of $L^2(X,\mu)$
orthogonal to 
constant functions, and by $\rho_0$ the regular representation of $G$
on $L^2_0(X,\mu)$. Now, with some abuse of terminology,  say that a sequence
$\{A_{t}
\mid {t}\in \bn\}$  of nontrivial measurable subsets of $X$ is  {\sl asymptotically
invariant\/}  if the
sequence of 
functions $1_{A_{t}} - \mu(A_{t})$  is {asymptotically
$\rho_0$-invariant}. Equivalently, if 
$$
{\mu(A_{t}\triangle
gA_{t})}/{\mu(A_{t})}\to 0\text{  as }{t}\to \infty \text{ uniformly on
compact subsets of }G\,.\tag AI
$$ 
Further,  we will say that $\{A_{t}\}$ is a {\sl $0$-sequence\/}
if  $\lim_{{t}\to \infty}\mu(A_{t}) = 0$. 

Now we can state the following useful criterion for  $\rho_0$ being
close to  $I_G$:

\proclaim{Proposition}  Let $G$ be a locally compact second countable
group acting ergodically on a probability space  $(X,\mu)$. Then the
following two conditions are equivalent: 

\roster
\item"{(i)}" there exists an asymptotically invariant $0$-sequence of
subsets $A_t$ of $X$; 

\item"{(ii)}" $\rho_0$ is not isolated from $I_G$.
\endroster 
\endproclaim

The implication (i)$\Rightarrow$(ii) is clear: by definition,  the sequence of
functions $1_{A_{t}} - \mu(A_{t})$ is  asymptotically 
$\rho_0$-invariant whenever $\{A_{t}\}$ is asymptotically invariant.
 K.~Schmidt \cite{S}, using a result of J.~Rosenblatt
\cite{Ro}, proved the converse for countable groups $G$; in fact, he
showed that both conditions are equivalent to

\roster
\item"{(iii)}"  $G$ has more than one invariant mean on $L^\infty(X,\mu)$.
\endroster 

In \cite{FS}, A.~Furman and Y.~Shalom extended the approach of
Rosenblatt and Schmidt to uncountable  groups. In particular, assuming
$G$ is locally 
compact, they
proved the implication (ii)$\Rightarrow$(iii), of which the converse
is  in this generality not always true. Our proof of Proposition A.0
 is based on the ideas  of
Rosenblatt-Schmidt-Furman-Shalom. However we have chosen to completely
avoid any use of invariant means, in the hope to make the argument
more transparent and less involved. 

\bigskip

\noindent{\it Proof of
Proposition A.0.\/}\  \ 
 Suppose we are given a sequence of functions
$\{\ph_{t}\}\in L^2_0(X,\mu)$ which is  asymptotically
$\rho_0$-invariant. Without loss of generality we can assume that all
the functions $\ph_{t}$ have $L^2$-norm $1$. Note also that any weak
limit point of the sequence $\{\ph_{t}\}$ must be $\rho_0$-invariant,
hence (by the ergodicity of the $G$-action on $X$) equal to zero. Thus, by
choosing a subsequence,  we can assume that $\ph_{t}\to 0$ weakly as
${t}\to\infty$.

Our goal is to produce
an asymptotically invariant $0$-sequence  $\{A_t\}$ of subsets of
$X$. 
Define a sequence $\{\sigma_{t}\}$ of
probability measures on $\br$ by
$$
\sigma_{t}(A) = \mu\big(\ph^{-1}(A)\big)\,,\quad A\subset \br\,.
$$
Observe that 
$$
\int_\br {z}\,d\sigma_{t}({z}) = 0\text{ and  }\int_\br
{z}^2\,d\sigma_{t}({z}) = 1\,.\tag A.0
$$
 In view of the last equality, 
we may
assume that $\sigma_{t}$ converges weakly on compacta to a probability
measure $\sigma$ on $\br$. The construction of the desired   sequence of 
sets will crucially depend on this measure. Following \cite{S} and
\cite{FS}, we  consider 
two cases.

\specialhead{Case 1.} {\sl The limit measure is concentrated on one
point ${a}\in\br$.}
\endspecialhead

\subhead{A.1.1}\endsubhead Let us, following \cite{FS}, first show that
${a} = 0$. Indeed, using (A.0), for any ${t}\in \bn$ and $N > 0$ one can write
$$
\left|\int_{-N}^N {z}\,d\sigma_{t}({z})\right| = \left|\int_{|{z}|>N}
{z}\,d\sigma_{t}({z})\right| = \frac1N\left|N\int_{|{z}|>N}
{z}\,d\sigma_{t}({z})\right| \le \frac1N\left|\int
{z}^2\,d\sigma_{t}({z})\right| = \frac1N\,.
$$ 
Choosing $N$ large enough and $\sigma_{t}$ close enough to $\sigma$, one
deduces that $|{a}| = \left|\int_{-N}^N {z}\,d\sigma({z})\right|$ must
be very small, which is only possible if ${a} = 0$. In particular,
this implies that for any $C > 0$,
$$
\int\limits_{\{|\ph_{t}| <
C\}} \ph_{t}^2\,d\mu = \int_{-C}^C {z}^2\,d\sigma_{t}({z}) \to \int_{-C}^C
{z}^2\,d\sigma({z}) = 0\,.\tag A.1.1
$$

\subhead{A.1.2}\endsubhead The next step is to pass from functions  
$\{\ph_{t}\}$ with zero mean value to another sequence $\{h_{t}\}$ of
nonnegative integrable functions. Namely we  define 
$$
h_{t}(x) = \cases \ph_{t}^2(x),\ &|\ph_{t}(x)|\ge 1\\ 0,\quad  &|\ph_{t}(x)| < 1
\endcases\tag A.1.2
$$

In what follows, we denote by $\|h\|_{1}$ the  $L^1$-norm of a function
$h$, and keep the notation $\|\cdot\|$ for the  $L^2$-norm.

\proclaim{Lemma} As ${t}\to\infty$, $\|h_{t}\|_{1}\to 1$ and
$\|h_{t}-gh_{t}\|_{1}\to 0$ uniformly on compact subsets of $G$.
\endproclaim

\demo{Proof} Note first that $\|\ph_{t}^2\|_{1} = \|\ph_{t}\|^2 = 1$, 
while $\|h_{t}\|_{1} - \|\ph_{t}^2\|_{1} = \int_{\{|\ph_{t}(x)| <
1\}}\ph_{t}^2\,d\mu \to 0$ in view of (A.1.1). Now for any $g\in G$ one
can write
$$
\|h_{t}-gh_{t}\|_{1} = \int\limits_{\{|g\ph_{t}| <
1,\,|\ph_{t}| \ge
1\}}\ph_{t}^2\,d\mu + \int\limits_{\{|\ph_{t}| <
1,\,|g\ph_{t}| \ge
1\}}g\ph_{t}^2\,d\mu + \int\limits_{\{|\ph_{t}| \ge
1,\,|g\ph_{t}| \ge
1\}}|\ph_{t}^2-g\ph_{t}^2|\,d\mu\,.
$$
The first integral in the r.h.s.~is not greater than
$$
 \dsize\int\limits_{\{1 \le |\ph_{t}| < 2
\}}\ph_{t}^2\,d\mu \ +\   \int\limits_{\{|\ph_{t}| \ge
2,\,|g\ph_t| \le |\ph_{t}|/2\}}\ph_{t}^2\,d\mu \ \le \int\limits_{\{1
\le |\ph_{t}| < 2 
\}}\ph_{t}^2\,d\mu \ +\   \frac43\int\limits_{\{|\ph_{t}| \ge
2\}}|\ph_{t}^2-g\ph_{t}^2|\,d\mu\,;
$$
 similarly, $\dsize\int\limits_{\{|\ph_{t}| <
1,\,|g\ph_{t}| \ge
1\}}g\ph_{t}^2\,d\mu \le \int\limits_{\{1 \le |g\ph_{t}| < 2
\}}g\ph_{t}^2\,d\mu \ +\   \frac43\int\limits_{\{|g\ph_{t}| \ge
2\}}|\ph_{t}^2-g\ph_{t}^2|\,d\mu$. Thus, using (A.1.1) and the
$G$-invariance of $\mu$, one gets
$$
\limsup_{{t}\to\infty}\|h_{t}-gh_{t}\|_{1} \le \frac{11}3\cdot
\limsup_{{t}\to\infty}\|\ph_{t}^2-g\ph_{t}^2\|_{1}\,.
$$
But $\|\ph_{t}^2-g\ph_{t}^2\|_{1} = \|(\ph_{t}-g\ph_{t})(\ph_{t}+g\ph_{t})\|_{1}
\le 2 \|\ph_{t}-g\ph_{t}\|$, and the latter $L^2$-norms tend to zero
uniformly on compact subsets of $G$, hence the claim.
\qed\enddemo

\subhead{A.1.3}\endsubhead The next step of the proof is to pass from functions to sets. Here we use the
following trick, dating 
back to I.~Namioka \cite{N}: if $h$ is a nonnegative
function on $X$ and $z \ge  0$, denote by $B_{z,h}$ the subset of $X$
given by 
$$
B_{z,h} \df \{x\in X\mid h(x) \ge z\}\,.
$$
Then one can  reconstruct the value of $h(x)$ as the Lebesgue measure
of the set $\{z\ge 0 \mid x\in B_{z,h}\}$. Moreover, if $g\in G$,
the absolute value of $(gh)(x) - h(x)$ is equal to the measure of
$\{z\ge 0 \mid x\in B_{z,h}\triangle B_{z,gh}\}$. Therefore,
assuming $h$ is integrable, its $L^1$-norm is equal to
$$
\|h\|_{1} = \int_X \int_0^\infty 1_{\{z \mid x\in B_{z,h}\}}\,dz\,d\mu(x)
= \int_0^\infty\int_X1_{\{z \mid x\in B_{z,h}\}}\,d\mu(x)\,dz = \int_0^\infty
\mu(B_{z,h}) \,dz \,;
$$
similarly,
$$
\|gh-h\|_{1} = \int_X \int_ 0^\infty1_{\{z \mid x\in B_{z,h}\triangle B_{z,gh}\}}\,dz\,d\mu(x)
 = \int_0^\infty
\mu(B_{z,h}\triangle B_{z,gh}) \,dz \,.
$$

This way, with $h_{t}$ as defined in (A.1.2), one deduces from Lemma A.1.2
that as ${t}\to\infty$,  
$$
 \int_0^\infty
\mu(B_{z,h_{t}}) \,dz \to 1 \text{ and }\int_0^\infty
\mu(B_{z,h_{t}}\triangle B_{z,gh_{t}}) \,dz \to 0\text{ uniformly on
compacta.}\tag A.1.3
$$
Furthermore, uniformly for all $z > 0$ one has 
$$
\mu(B_{z,h_{t}}) = \mu\big(\{x\mid h_{t}(x) \ge z\}\big) \le \mu\big(\{
x\bigm| |\ph_{t}(x)| \ge 1\}\big) = \sigma_{t}\big(\br\ssm(-1,1)\big) \to
0\,,\tag A.1.4
$$
since by assumption the limit measure is concentrated at $0$.

\subhead{A.1.4}\endsubhead The final step is to get rid of integration
over $z$ in (A.1.3). Choose a 
sequence $\{K_{t}\mid {t}\in\bn\}$ of compact subsets of $G$ such that:

\roster
\item"{(i)}" $e\in K_{t}$ for all ${t}$;

\item"{(ii)}" $K_{t} \subset K_{{t}+1}$ for all ${t}$, and $\cup_{{t} = 1}^\infty K_{t} = G\,$;

\item"{(iii)}" each $K_{t}$ is equal to the closure of its interior.
\endroster

\noindent Fix  a right-invariant Haar measure $\nu$  on $G$. From (i) and (iii)
it follows that for any 
${t}$ the value of $\inf_{g\in
K_{t}}\frac{\nu(K_{t}\cap K_{t}g)}{\nu(K_{t})}$ is positive. Thus one can
choose a sequence of positive numbers $\vre_{t}$ with $\vre_{t} \to 0 $ as
${t}\to\infty$ such that 
$$
\nu(K_{t}\cap K_{t}g) \ge \vre_{t} \nu(K_{t})\text{ for all }g\in
K_{t}\,.\tag A.1.5
$$
 Now,  replacing $\{h_{t}\}$ by a subsequence if
needed, in view of (A.1.3) we can assume that for all $g\in
K_{t}$
$$
 \int_0^\infty
\mu(B_{z,h_{t}}\triangle B_{z,gh_{t}}) \,dz  < \frac{\vre_{t}^2}4 \int_0^\infty
\mu(B_{z,h_{t}}) \,dz\,.
$$
Integrating over $K_{t}$ and then changing the order of integration
between $dz$ and $d\nu$, we find that
$$
 \int_0^\infty\int_{K_{t}} \left( \tfrac{\vre_{t}^2}4\mu(B_{z,h_{t}}) -
\mu(B_{z,h_{t}}\triangle B_{z,gh_{t}}) \right)\,d\nu(g)\,dz > 0\,.
$$
Therefore for every ${t}$ there exists $z_{t} > 0$ such that
$$
\frac1{\nu(K_{t})}\int_{K_{t}}\mu(B_{z_{t},h_{t}}\triangle B_{z_{t},gh_{t}})\,d\nu(g)
<  \frac{\vre_{t}^2}4 \mu(B_{z_{t},h_{t}})\,.\tag A.1.6
$$
Let us now show that the sets $A_{t} \df
B_{z_{t},h_{t}}$ form an asymptotically invariant $0$-sequence.
It is immediate from (A.1.4) that $\mu(B_{z_{t},h_{t}}) \to 0$ as
${t}\to\infty$. Thus it suffices to   find a sequence of compacta
$\{K'_{t}\}$ exhausting  $G$ 
such that  
$$
\mu(A_{t}\triangle gA_{t}) / \mu(A_{t}) \le \vre_{t}\text{ whenever }
g\in K'_{t}\,.\tag A.1.7
$$
 This is achieved by putting $K'_{t} \df  Q_{t}^{-1} Q_{t}$, where
$$
Q_{t} \df \{g\in K_{t}\mid \mu(A_{t}\triangle gA_{t}) \le \frac{\vre_{t}}2
\mu(A_{t})\}\,.\tag A.1.8
$$
(Indeed, if $g = g_1^{-1}g_2$, with  $g_1,g_2\in
Q_{t}$, then $\mu(A_{t}\triangle gA_{t}) = \mu(g_1A_{t}\triangle g_2A_{t}) \le
\mu(A_{t}\triangle g_1A_{t}) + \mu(A_{t}\triangle g_2A_{t})$, and
(A.1.7) follows.) Therefore, the claim for Case 
1 can be derived  from condition (ii) and  the following 
 
\proclaim{Lemma} $K'_{t}$ contains $K_{t}$.
\endproclaim

\demo{Proof}  If not, then there exists $g\in K_{t}$ such that $Q_{t}g\cap
Q_{t} = \varnothing$, which implies that $Q_{t}g \subset (K_{t} \ssm Q_{t})
\cup (K_{t} g \ssm K_{t})$ $\thus$ $\nu(Q_{t}) \le \nu(K_{t}) - \nu(Q_{t})
+\nu(K_{t}) -  
\nu(K_{t}\cap K_{t}g)$. This, in view of (A.1.5), forces $\nu(Q_{t})$ to be
not greater than $(1 - \frac{\vre_{t}}2)\nu(K_{t})$. On the other hand,
using (A.1.8) and (A.1.6), one can write 
$$
\frac{\vre_{t}}2 \mu(A_{t})\nu(K_{t} \ssm Q_{t}) < \int_{K_{t}\ssm Q_{t}}\mu(A_{t}\triangle
gA_{t})\,d\nu(g) \le  \int_{K_{t}}\mu(A_{t}\triangle 
gA_{t})\,d\nu(g) < \frac{\vre_{t}^2}4 \mu(A_{t})\nu(K_{t})\,,
$$
therefore $\nu(K_{t} \ssm Q_{t}) < \frac{\vre_{t}}2 \nu(K_{t})$, a
contradiction. \qed\enddemo

\specialhead 
{Case 2.} {\sl The limit measure $\sigma$ is not
concentrated on one 
point.}
\endspecialhead

\subhead{A.2.1}\endsubhead The above assumption implies that there
exists ${a} \in\br$ such that
$$
0 < \sigma\big(({a},\infty)\big) = \sigma\big([{a},\infty)\big)
= \tau < 1\,.\tag A.2.1
$$
Without loss of generality we can assume that ${a} > 0$. As a first
attempt to build a good sequence of sets out of $\{\ph_{t}\}$, we
consider $B_{t}\df 
\ph_{t}^{-1}\big(({a},\infty)\big)$. Then clearly  $\mu(B_{t})\to\tau$
as ${t}\to\infty$. Moreover, one has

\proclaim{Lemma} The sequence $\{B_{t}\}$ is asymptotically invariant.
\endproclaim

\demo{Proof} In view of  (A.2.1), for any $\vre > 0$ one can find 
$\delta > 0$ such that $
\sigma\big(({a}-\delta,{a}+\delta)\big) = $ 
$\mu\big(\big\{x\bigm||\ph_{t}(x) - {a}| < \delta\big\}\big) \le
\vre$. Then $\mu(B_{t}\triangle gB_{t})$ is not greater than
$$
\mu\big(\big\{x\bigm||\ph_{t}(x)
- {a}| < \delta\big\}\big) + \mu\big(\big\{x\bigm||\ph_{t}(x) -
{a}| \ge \delta\big\} \cap (B_{t}\triangle gB_{t}) \big) \le \vre +
\frac1{\delta^2}\int_X|g\ph_{t} - \ph_{t}|^2\,d\mu\,. 
$$
Since $\{\ph_{t}\}$  is
asymptotically 
$\rho_0$-invariant, 
$\limsup_{{t}\to\infty}\mu(B_{t}\triangle gB_{t}) \le \vre$ uniformly on compacta, and (AI) follows.
 \qed\enddemo

\subhead{A.2.2}\endsubhead We now use $\{B_{t}\}$  to produce  a
family of  asymptotically invariant sequences 
$
B_{t}^{(k)}$ with $\limsup_{{t}\to \infty}\mu(B_{t}^{(k)}) \to 0$ as
$k\to\infty$. As a first step, choose a  
sequence $l_{t}\to\infty$ and  a  
sequence of increasing compact subsets $K_{t}$ of $G$ exhausting $G$
such that

\roster
\item"{(i)}"
 $\mu(B_{l_{t}})  \le \tau + 1/t$;
\item"{(ii)}" $\mu(B_{l}\triangle gB_{l}) \le 1/2t$ uniformly  in $g\in
K_{t}$ for all $l \ge l_{t}$;
\item"{(iii)}" $\#\{l\mid \mu(B_{l_{t}}\cap B_{l}) > 0\} = \infty$.
\endroster

 \noindent Observe that by the Schwarz inequality, for any $l$ one has
$$
\left|\int_{B_{l_{t}}\ssm B_l} \ph_l\,d\mu\right| \le
\left(\int_{B_{l_{t}}\ssm B_l}
\ph_l^2\,d\mu\right)^{1/2}\left(\int_{B_{l_{t}}\ssm B_l}
1\,d\mu\right)^{1/2} \le \sqrt{\mu(B_{l_{t}}\ssm B_{l})}\,.
$$
Therefore
$$
\left|\int_{B_{l_{t}}} \ph_l\,d\mu\right| 
\ge
\int_{B_{l_{t}}\cap B_l} \ph_l\,d\mu - \left|\int_{B_{l_{t}}\ssm B_l}
\ph_l\,d\mu\right| \ge {a}  \mu(B_{l_{t}}\cap
B_{l})-\sqrt{\mu(B_{l_{t}}\ssm B_{l})}\,.\tag A.2.2
$$
Applying (iii) and the weak convergence of $\{\ph_{t}\}$ to zero,
for each ${t}$ choose  $l > l_{t}$ such that $z \df \mu(B_{l_{t}}\cap
B_{l}) > 0$ and 
$
\left|\int_{B_{l_{t}}} \ph_l\,d\mu\right| <
\frac{a}2\mu(B_{l_{t}})\,.$
Combining this with (A.2.2), we obtain the 
inequality ${a} z -
\sqrt{\mu(B_{l_{t}}) - z} < \frac{a}2\mu(B_{l_{t}})$. An exercise in
quadratic equations gives that $z$ must be less than
$\frac{\mu(B_{l_{t}})}2 + \frac{\sqrt{1 + 2{a}^2\mu(B_{l_{t}})} -
1}{2{a}^2}$. 

Now denote $B_{t}^{(2)}\df B_{l_{t}}\cap B_l$. Then $\limsup_{{t}\to
\infty}\mu(B_{t}^{(2)})\le \tau^{(2)} \df \frac{\tau}2 + \frac{\sqrt{1
+ 2{a}^2\tau} - 
1}{2{a}^2}$. Also, from (ii) it follows that $\mu(B_{l}\triangle
gB_{l}) \le 1/2t$ and $\mu(B_{l_{t}}\triangle
gB_{l_{t}}) \le 1/2t$ uniformly  in $g\in
K_{t}$. Therefore $\mu(B_{t}^{(2)}\triangle
gB_{t}^{(2)}) \le \mu(B_{l}\triangle
gB_{l}) + \mu(B_{l_{t}}\triangle
gB_{l_{t}}) \le 1/t$, which shows that $\{B_{t}^{(2)}\}$ is asymptotically invariant.

Applying the above procedure to $\{B_{t}^{(2)}\}$ we produce another
sequence $B_{t}^{(3)}\df B_{l_{t}}^{(2)}\cap B_l$ for appropriate ${l_{t}}$
and 
$l > {l_{t}}$, and, inductively, a family of  asymptotically invariant sequences
$
B_{t}^{(k)}\df
B_{l_{t}}^{(k-1)}\cap B_l\,, 
\quad\text{with}\quad\limsup_{{t}\to
\infty}\mu(B_{t}^{(k)})\le \tau^{(k)} \df \frac{\tau^{(k-1)}}2 + \frac{\sqrt{1 +
2{a}^2\tau^{(k-1)}} - 
1}{2{a}^2}\,.
$
 It is easy to see that $\tau^{(k)}\to 0$ as
$k\to\infty$. Finally,  define $A_{t}$ diagonally as
$B_{{t}'}^{({t})}$, where 
${{t}' > {t}}$ is chosen so that $\mu(A_{t}\triangle
gA_{t}) < 1/{t}$ whenever $g$ belongs to the compact set $K_{t}$. This
completes the construction of the asymptotically invariant $0$-sequence
$\{A_{t}\}$,  as well
as the proof of Proposition A.0.  \qed

\subhead{A.3}\endsubhead It remains to write down the 

\demo{Proof of Lemma 3.1} It is easy to deduce  from (3.1) and the
$G$-equivariance 
of $\pi$ that if  $\{A_{t}\}$ is
an asymptotically invariant $0$-sequence of subsets of $X_1$, then
$\{\pi(A_{t})\}$ is 
an asymptotically invariant $0$-sequence of subsets of $X_2$; and,
conversely, any asymptotically invariant $0$-sequence 
$\{ A_{t}\}$ of subsets of $X_2$ gives rise to an asymptotically
invariant $0$-sequence 
$\{\pi^{-1}(A_{t})\}$ of subsets of $X_1$. \qed\enddemo

\pagebreak

%
 

\heading{  Acknowledgements}\endheading

The authors want to thank Marc Burger, Nikolai Chernov, Dmitry Dolgopyat, Yves
Guivarc'h, Alex Furman,
David Kazhdan, Yuval Peres, Marc Pollicott, Peter
Sarnak and Yehuda
Shalom  for helpful discussions, and the referee for useful
comments. Thanks are also due to organizers 
and participants of the Workshop on Geometric Rigidity and Hyperbolic Dynamics
(Oberwolfach, June 1996) where some of the results of this paper were
announced.


\vskip .3in


\Refs
\widestnumber\key{Bou2}


\ref\key B \by M. Bekka \paper On uniqueness of invariant means \jour
Proc. Amer. Math. Soc. \vol 126 \pages 507--514 \yr 1998  \endref

\ref\key Bo \by A. Borel \book Introduction aux groupes
arithm\'etiques  \publ Hermann \publaddr Paris \yr 1969 \endref

\ref\key Bos \by M. Boshernitzan \paper Quantitative recurrence
results \jour Invent. Math. \vol 113 \yr 1993 \pages 617--631 \endref

\ref\key Bou1 \by N. Bourbaki \book  Integration, Chapitres VII et VIII \publ Hermann
\publaddr Paris \yr 1963 \endref

\ref\key Bou2 \bysame \book 
Groupes et Alg\`ebres de Lie, Chapitres IV, V et VI\publ Hermann
\publaddr Paris \yr 1968 \endref



\ref\key BS  \by M. Burger and P. Sarnak\paper Ramanujan Duals II \jour
Inv. Math.\vol 106\pages 1--11\yr 1991\endref

\ref\key CK \by N. Chernov and D. Kleinbock \paper Dynamical
Borel-Cantelli lemmas for Gibbs measures \paperinfo Preprint \yr 1999
\endref 
 
\ref\key CR \by J.-P. Conze and A. Raugi \paper Convergence des
potentiels pour un op\'erateur de transfert, 
applications aux syst\`emes dynamiques et aux cha\^\i nes de
Markov\paperinfo Preprint \yr 1999 
\endref 
 
\ref\key D  \by S. G. Dani \paper Divergent trajectories of flows on
\hs s and \da \jour
J. Reine Angew. Math.\vol 359\pages 55--89\yr 1985\endref

\ref\key Do  \by M. M. Dodson \paper Geometric and probabilistic ideas
in metric \da\jour
Russian Math. Surveys\vol 48\pages 73--102\yr 1993\endref

\ref\key Dol  \by D. Dolgopyat 
%
\paper Limit theorems for partially hyperbolic
systems\paperinfo Preprint \yr 1999 
\endref 




\ref\key F  \by Y. Flicker \paper Automorphic forms on covering groups
of $GL(2)$ \jour
Inv. Math.\vol 57\pages 119--182\yr 1980\endref

\ref\key FS \by A. Furman and Y. Shalom \paper Sharp ergodic theorems
for group actions and strong ergodicity \jour
Erg. Th. Dyn. Systems \toappear \endref

\ref\key G \by A. V. Groshev  \paper Une th\'eor\`eme sur les syst\`emes des formes lin\'eaires \jour Dokl. Akad. Nauk SSSR \vol 9 \yr 1938 \pages 151--152 \endref 
 
\ref \key GR \by H. Garland and  M. S. Raghunathan \paper Fundamental
domains for lattices in R-rank $1$ semisimple Lie groups \jour
Ann. Math. \vol 92 \pages 279--326 \yr 1970 \endref

\ref\key HV \by R. Hill and S. Velani \paper Ergodic theory of
shrinking targets \jour  Invent. Math. \vol 119 \yr 1995 \pages
175--198 \endref 


\ref\key J \by L. Ji\paper Metric compactifications of locally
symmetric spaces \jour Internat. J. Math. \vol 9 \yr 1998 \pages
465--491\endref 

\ref\key KM \by D.  Kleinbock and G. A. Margulis \paper  Bounded
orbits of nonquasiunipotent flows  on homogeneous spaces \jour
Amer. Math. Soc. Transl. \vol 171 \pages 141--172 \yr 1996 
 \endref 



\ref\key KS \by S. Kochen and C. Stone \paper A note on the Borel Cantelli lemma \jour Ill. J. Math. \vol 8 \pages 248--251 \yr 1964 \endref

\ref\key L \by E. Leuzinger\paper Tits geometry and arithmetic groups
\paperinfo Preprint \yr 1997 \endref 


\ref\key M \by G. A. Margulis \book  Discrete subgroups of semisimple Lie groups
 \publ Springer-Verlag \publaddr Berlin and New York \yr 1991 \endref  


\ref\key Ma\by F. J. Mautner \paper Geodesic flows on symmetric Riemannian spaces \jour Ann. Math. \vol 65\pages 416--431\yr 1957\endref%
 
\ref\key N \by I. Namioka \paper Folner's conditions for amenable
semi-groups \jour Math. Scand. \vol 15\yr 1964 \pages 18--28\endref

\ref\key P \by W. Philipp \paper Some metrical theorems in number theory \jour Pacific J. Math. \vol 20 \pages 109--127 \yr 1967 \endref


\ref\key R1 \by M. S. Raghunathan \book Discrete subgroups of Lie groups
\publ Springer-Verlag \publaddr Berlin and New York \yr 1972 \endref


\ref\key R2  \bysame \paper On the
congruence  subgroup problem\jour
Publ. Math. IHES\vol 46\pages 107--161\yr 1976\endref


\ref\key Ro  \by J. Rosenblatt \paper Uniqueness of invariant means
for measure preserving transformations \jour
Trans. Amer. Math. Soc. \vol 265\pages 623--636\yr 1981\endref

\ref\key S  \by K. Schmidt \paper Amenability, Kazhdan's property
(T), strong ergodicity and invariant means for ergodic group actions
\jour Erg. Th. Dyn. Systems 
\vol 1\pages 223--236\yr 1981\endref


\ref\key Si \by C. L. Siegel \paper A mean value theorem in geometry
of numbers \jour Ann. Math. \vol 46 \pages 340--347 \yr 1945 \endref


\ref\key Sk \by M. Skriganov \paper Ergodic theory on $SL(n)$, \da s
and anomalies in the lattice point problem  \jour  Invent. Math. \vol
132 \yr 1998 \pages 1--72 \endref 

\ref\key Sp \by F. Spitzer \book Principles of random walk \publ Van Nostrand \publaddr Princeton \yr 1964 \endref


\ref\key {Spr}\by V. \spr \book Metric theory of Diophantine
approximations \publ
John Wiley \& Sons \publaddr New York-Toronto-London \yr 1979\endref



\ref\key Su \by D. Sullivan \paper Disjoint spheres, approximation by imaginary quadratic numbers, and the logarithm law for geodesics \jour Acta Math. \pages 215--237 \vol 149 \yr 1982 \endref

\ref\key T1 \by J. Tits\paper Algebraic and abstract simple groups
\jour Ann. Math. \vol 80 \yr 1964 \pages 
313--329\endref 

\ref\key T2 \bysame \paper Classification of algebraic semisimple
groups \inbook 1966 Algebraic Groups
and Discontinuous Subgroups (Proc. Sympos. Pure Math., Boulder, Colo.,
1965) \pages 33--62 \publ Amer.
Math. Soc. \publaddr Providence, R.I. \yr 1966 \endref

\ref\key V \by M.-F. Vigneras \paper Quelques remarques sur la
conjecture $\lambda_1 \ge \frac14$ \inbook Seminar in number theory,
Paris 1981-82, Progr. Math., 38 \publ Birkh\" auser Boston \publaddr
Boston, MA \yr 1983 \pages 321--343 \endref

\ref\key W \by G. Warner \book Harmonic analysis on semisimple Lie groups I \publ Springer-Verlag \publaddr Berlin and New York \yr 1972 \endref




\ref\key Z \by R. Zimmer \book  Ergodic theory and
semisimple groups
 \publ Birkh\"auser \publaddr Boston  \yr 1984 \endref  

\endRefs
\enddocument